\documentclass[12pt]{amsart}

\usepackage{amsmath}
\usepackage{amsfonts}
\usepackage{amssymb}
\usepackage[all]{xy}           %xypic macro for latex2.09

\usepackage{txfonts}
\usepackage{amscd}

\usepackage[shortlabels]{enumitem}

\usepackage{tikz}

\makeatletter

\newtheorem{thm}{Theorem}[section]
\newtheorem{lem}[thm]{Lemma}
\newtheorem{cor}[thm]{Corollary}
\newtheorem{pro}[thm]{Proposition}
\newtheorem{ex}[thm]{Example}
\newtheorem{rmk}[thm]{Remark}
\newtheorem{defi}[thm]{Definition}

\newcommand {\emptycomment}[1]{}

\newcommand{\hd}{\mathcal{D}}
\newcommand{\g}{\mathfrak g}

\newcommand{\gl}{\mathfrak {gl}}

\newcommand{\cb}{\mathbb C}

\newcommand{\fd}{\mathbf d}

\newcommand{\Der}{\mathrm{Der}}

\newcommand{\Inn}{\mathrm{Inn}}

\newcommand{\Aut}{\mathrm{Aut}}

\newcommand{\Ext}{\mathrm{Ext}}
\newcommand{\ad}{\mathrm{ad}}

\def\oll#1#2{\mathrel {{\triangleright}_{(#1)}^{#2}}}
\def\orr#1#2{\mathrel {{\triangleleft}_{(#1)}^{#2}}}
\def\obl#1{\mathrel {{\blacktriangleright}_{(#1)}}}
\def\bll#1#2{\mathrel {{\blacktriangleright}_{(#1)}^{#2}}}
\def\obll#1{\mathrel {{\bar{\blacktriangleright}}_{(#1)}}}
\def\ool#1{\mathrel {{\triangleright}_{(#1)}}}
\def\obr#1{\mathrel {{\blacktriangleleft}_{(#1)}}}
\def\brr#1#2{\mathrel {{\blacktriangleleft}_{(#1)}^{#2}}}
\def\obrr#1{\mathrel {{\bar{\blacktriangleleft}}_{(#1)}}}
\def\oor#1{\mathrel {{\triangleleft}_{(#1)}}}
\def\oo#1{\mathrel {{}_{(#1)}}}
\def\ooo#1#2{\mathrel {{}_{(#1)}^{#2}}}

\def\id{\mathop {\fam0 id}\nolimits}

\def\Ker{\mathrm{Ker}\,}
\def\Cok{\mathrm{Coker}\,}
\def\Img{\mathrm{Im}\,}
\def\Hom{\mathrm{Hom}}

\def\HH{\mathrm{HH}}

\begin{document}

\title[Crossed modules, non-abelian extensions and Wells exact sequence]
{Crossed modules, non-abelian extensions of associative conformal algebras and Wells exact sequences}

\author{Bo Hou}
\address{School of Mathematics and Statistics, Henan University, Kaifeng 475004,
China}
\email{bohou1981@163.com}
\vspace{-5mm}

\author{Jun Zhao}
\address{School of Mathematics and Statistics, Henan University, Kaifeng 475004,
China}
\email{zhaoj@henu.edu.cn}

%\date{\today}

\begin{abstract}
In this paper, we introduce the notions of crossed module of associative conformal algebras,
2-term strongly homotopy associative conformal algebras, and discuss the relationship
between them and the 3-th Hochschild cohomology of associative conformal algebras.
We classify the non-abelian extensions by introducing the non-abelian
cohomology. We show that non-abelian extensions of an associative conformal algebra
can be viewed as Maurer-Cartan elements of a suitable differential graded Lie algebra,
and prove that the Deligne groupoid of this differential graded Lie algebra corresponds
one to one with the non-abelian cohomology. Based on this classification,
we study the inducibility of a pair of automorphisms about a non-abelian extension of
associative conformal algebras, and give the fundamental sequence of Wells in the context of
associative conformal algebras. Finally, we consider the extensibility of a pair of
derivations about an abelian extension of associative conformal algebras,
and give an exact sequence of Wells type.
\end{abstract}


\maketitle

%\tableofcontents %目录

%\setcounter{section}{0}

%\allowdisplaybreaks

%\end{document}
\vspace{-4mm}
%%%%%%%%%%%%%%%%%%%%%%%%%%%%%%%%%%%%%%%%%%%%%%%%%%%%%%%%%%%%%%%%%%%%%%%%%%%%%%%%%
%    section  1   Introduction
%%%%%%%%%%%%%%%%%%%%%%%%%%%%%%%%%%%%%%%%%%%%%%%%%%%%%%%%%%%%%%%%%%%%%%%%%%%%%%%%%%%%%%

\section{Introduction}\label{sec:intr}

The notion of a conformal algebra encodes an axiomatic description
of the operator product expansion of chiral fields in conformal field
theory. The theory of Lie conformal algebras appeared as a formal language
describing the algebraic properties of the operator product
expansion in two-dimensional conformal field theory (\cite{BPZ,Ka,Ka1}).
The structure of a Lie conformal algebra gives an axiomatic description of the
operator product expansion of chiral fields in conformal field theory.
Associative conformal algebras naturally come from representations
of Lie conformal algebras. Moreover, some Lie conformal algebras appeared
in physics are embeddable into an associative one \cite{Ro, Ro1}.
The structure theory and representation theory of associative conformal algebras have
attracted the attention of many scholars and achieved a series of results, see
\cite{BFK1,BFK2,Ko,Ko1,Re,Re1}.

The concept of crossed modules was introduced by J.H.C. Whitehead in the late 1940s \cite{Whi}.
Crossed modules of Lie algebras first appeared in the work of Gerstenhaber \cite{Ge1}.
Later, the crossed modules of various algebraic structures have been extensively
studied \cite{BM,CLP,CKL1,CG,LL,Wag}. They have therefore been used in homotopy theory,
homological algebra, non-abelian cohomology, algebraic K-theory, ring theory, combinatorial
group theory and applications of the related algebra \cite{BH,Hol,Lod}.
Crossed modules of groups and Lie algebras turn up in the book of homological algebra
as interpretations of cohomology classes of cohomological degree 3. It is well-known that
the category of crossed modules of Lie algebras is equivalent to the category of strict
Lie 2-algebras \cite{BC}. Here we introduce crossed modules of associative
conformal algebras and consider conformal analogue of these results.

Extension problem is a famous and still open problem. The non-abelian extension theory is a
relatively general one in various kinds of extensions. In \cite{EM}, Eilenberg and Maclane
first considered non-abelian extensions of abstract groups. Subsequently, the non-abelian
extension theory has been widely studied in various fields of mathematics,
to which a vast literature was devoted, such as Lie groups \cite{Nee}, Leibniz algebras
\cite{CKL,LSW}, Lie algebras \cite{FP,Fr,IKL}, associative algebras \cite{AM,Gou},
3-Lie algebras \cite{SMT}, Lie 2-algebras \cite{TS}, Hom-Lie algebras \cite{ST,ST1},
Rota-Baxter algebras \cite{DHM,DR}, groupoids \cite{MOW}, etc.
Another interesting study related to extensions of algebraic
structures is given by the inducibility of pair of automorphisms.
Such study was first initiated by Wells in extensions of abstract groups \cite{Wel}.
In \cite{Wel}, the author defined a map, known as Wells map, and constructed a short
exact sequence (the fundamental short exact sequence of Wells) connecting various
automorphism groups. The Wells map and the inducibility of a pair of
automorphisms about an abelian extension were studied in the context of groups \cite{JL},
Lie (super)algebras \cite{BS,HH}, 3-Lie algebras \cite{TX}, Lie coalgebra \cite{DT}.
Recently, Das and his partners have studied the Wells map and the inducibility of a pair
of automorphisms about a non-abelian extension of Rota-Baxter algebras \cite{DHM,DR}.
Not only automorphisms, but extensibility of derivations about an abelian extension of Lie
algebras and associative algebras have been studied in \cite{BS1} and \cite{TX1}.
But little is known about the non-abelian extension of conformal algebras.
In this and the following paper \cite{ZH}, we will study the abelian extensions and the
non-abelian extensions of conformal algebras.

In this paper, we first introduce the notion of 2-term strongly homotopy associative
conformal algebras and obtain a 1-1 correspondence between equivalence
classes of skeletal 2-term strongly homotopy associative conformal algebras and
3-th Hochschild cohomologies of associative conformal algebras. We introduce the notion
of crossed modules of associative conformal algebras and give a 1-1 correspondence
between strict 2-term strongly homotopy associative conformal algebras and
crossed modules of associative conformal algebras. Moreover, we also relate the
3-th Hochschild cohomology of associative conformal algebras by means of a particular
kind of crossed modules.
Second, we consider non-abelian extensions of associative conformal algebras.
For given two associative conformal algebras $A$ and $B$, we define the non-abelian
2-cocycle on $B$ with values in $A$, and give an equivalence relation on the set
$\mathcal{Z}^{2}_{nab}(B, A)$ of all the non-abelian 2-cocycle on $B$ with values in $A$.
The non-abelian cohomology of $B$ with values in $A$, denoted by $\HH^{2}_{nab}(B, A)$,
is the quotient of $\mathcal{Z}^{2}_{nab}(B, A)$ by this equivalence relation.
By giving a 1-1 correspondence between the equivalence class of non-abelian extension
of $B$ by $A$ and the equivalence class of $\mathcal{Z}^{2}_{nab}(B, A)$,
we prove that the non-abelian extensions of $B$ by $A$ can be classified by
$\HH^{2}_{nab}(B, A)$. For this classification, we can also use the Maurer Cartan elements
of a suitable differential graded Lie algebra to characterize.
For the direct product conformal algebra $A\oplus B$, there is a differential graded
Lie algebra $(\mathcal{C}^{\bullet+1}(A\oplus B, A\oplus B),\; [-,-],\; d_{\ast})$ on
the Hochschild type complex of $A\oplus B$. We construct a
differential graded Lie subalgebra $\mathcal{L}$ of $(\mathcal{C}^{\bullet+1}(A\oplus B,
A\oplus B),\; [-,-],\; d_{\ast})$, and show that the Maurer Cartan elements of
$\mathcal{L}$ control the non-abelian extensions of $B$ by $A$.
Based on this results, we define the Wells map $\mathcal{W}:\; \Aut(A)\times\Aut(B)\rightarrow
\HH^2_{nab}(B, A)$, and give a necessary and sufficient condition for a pair of
automorphisms in $\Aut(A)\times\Aut(B)$ to be inducible by the Wells map. Moreover,
using the Wells map,  we obtain the Wells short exact sequence connecting of various
automorphism groups and cohomology in the context of associative conformal algebras.
Since the non-abelian extension is the generalization of abelian extension,
here we also give the corresponding results of abelian extension for each obtained
results of non-abelian extension. Moreover, for an abelian extension of associative
conformal algebras, we also define a Wells map, give a necessary and sufficient condition
for a pair of derivations to be extensible, and obtain an exact sequence of Wells type.

The paper is organized as follows.
In Section \ref{sec:prel}, we recall the notions of associative conformal
algebras, and the differential graded Lie algebra structure on the Hochschild complex
of associative conformal algebras. In Section \ref{sec:cross},
we introduce the notions of crossed module of associative conformal algebras
and 2-term strongly homotopy associative conformal algebras. We will discuss the
relationship between them and the 3-th Hochschild cohomology of associative conformal
algebras with coefficient in bimodules. In Section \ref{sec:Non-abel}, we define a
non-abelian cohomology group, and show that the non-abelian extensions can be
classified by the non-abelian cohomology group. In Section \ref{sec:Deli},
we construct a differential graded algebra $\mathcal{L}$ for two
associative conformal algebras $A$ and $B$. We show that the non-abelian
cocycles are in bijection with the Maurer-Cartan elements of $\mathcal{L}$,
and get that the equivalence relation on $\mathcal{Z}^{2}_{nab}(B, A)$ can be
interpreted as gauge equivalence relation on the set ${\rm MC}(\mathcal{L})$
of Maurer Cartan elements of $\mathcal{L}$.
In Section \ref{sec:ind-aut}, we study the inducibility of a pair of automorphisms
about a non-abelian extension of associative conformal algebras by the Wells map,
and obtain the fundamental sequence of Wells in the context of associative conformal algebras.
In Section \ref{sec:ind-der}, for an abelian extension of associative conformal algebras,
we consider the extensibility of a pair of derivations
about this extension, and give an exact sequence of Wells type.

Throughout this paper, we fix $\cb$ an algebraical closed field and characteristic zero
(for example, the field of complex numbers). All vector spaces are $\cb$-vector spaces,
all linear maps and bilinear maps are $\cb$-linear, all tensor product are over $\cb$,
unless otherwise specified. For any vector space $V$ and variable $\lambda$,
we use $V[\lambda]$ to denote the set of polynomials of $\lambda$ with coefficients in $V$.

%%%%%%%%%%%%%%%%%%%%%%%%%%%%%%%%%%%%%%%%%%%%%%%%%%%%%%%%%%%%%%%%%%%%%%%%%%%%%%%%%
%    section  2   Associative conformal algebras
%%%%%%%%%%%%%%%%%%%%%%%%%%%%%%%%%%%%%%%%%%%%%%%%%%%%%%%%%%%%%%%%%%%%%%%%%%%%%%%%%%%%%%
\section{Associative conformal algebras and Hochschild cohomology} \label{sec:prel}

We recall the notions of associative conformal algebras, conformal bimodule over
associative conformal algebras, and the Hochschild cohomology of an associative
conformal algebra with coefficients in a bimodule. For the details see
\cite{BDK,DK,BKV,HSZ}.

\begin{defi}\label{de:alg}
A {\rm conformal algebra} $A$ is a $\cb[\partial]$-module endowed with a
bilinear map $\cdot\oo{\lambda}\cdot : A\times A\rightarrow A[\lambda]$,
$(a, b)\mapsto a\oo{\lambda} b$ satisfying
\begin{eqnarray*}
\partial a\oo{\lambda}b=-\lambda a\oo{\lambda}b,  \qquad
a\oo{\lambda}\partial b=(\partial+\lambda)a\oo{\lambda}b,
\end{eqnarray*}
for any $a, b\in A$.
An {\rm associative conformal algebra} $A$ is a conformal algebra satisfying
\begin{eqnarray*}
(a\oo{\lambda}b)\oo{\lambda+\mu}c=a\oo{\lambda}(b\oo\mu c),
\end{eqnarray*}
for any $a, b, c\in A$.
\end{defi}

Let $(A,\; \cdot\ooo{\lambda}{A}\cdot)$, $(B,\; \cdot\ooo{\lambda}{B}\cdot)$ be two
associative conformal algebras. A $\cb[\partial]$-module homomorphism
$f: A\rightarrow B$ is call a {\it homomorphism of associative conformal algebras}
if for any $a, b\in A$, $f(a\ooo{\lambda}{A}b)=f(a)\ooo{\lambda}{B}f(b)$.
A homomorphism $f: A\rightarrow B$ is said to be an {\it isomorphism} if $f$ is a
bijection. For an associative conformal algebra $A$, we denote the automorphism goup
by $\Aut(A)$.

\begin{ex}
Let $(A,\cdot)$ be an associative algebra. Then ${\rm Cur}(A)=\cb[\partial]\otimes A$
is an associative conformal algebra with the following $\lambda$-product:
\begin{eqnarray*}
(p(\partial)a)\oo\lambda (q(\partial)b)=p(-\lambda)q(\lambda+\partial)(a\cdot b),
\end{eqnarray*}
for any $p(\partial), q(\partial)\in\cb[\partial]$ and $a, b\in A$.
\end{ex}

Now we recall that definition of left (or right) modules over an associative
conformal algebra.

\begin{defi}\label{def:mod}
A {\rm (conformal) left module} $M$ over an associative conformal algebra $A$ is a
$\cb[\partial]$-module endowed with a $\cb$-bilinear map $A\times M\rightarrow
M[\lambda]$, $(a, v)\mapsto a\ool{\lambda} v$, satisfying the following axioms:
\begin{align*}
(\partial a)\ool{\lambda} v&=-\lambda a\ool{\lambda} v,\qquad
a\ool{\lambda}(\partial v)=(\partial+\lambda)(a\ool{\lambda} v),\\
&(a\oo{\lambda} b)\ool{\lambda+\mu}v=a\ool{\lambda}(b\ool{\mu} v),
\end{align*}
for any $a, b\in A$ and $v\in M$. We denote it by $(M, \triangleright)$.

A {\rm (conformal) right module} $M$ over an associative conformal algebra $A$ is a
$\cb[\partial]$-module endowed with a $\cb$-bilinear map $M\times A\rightarrow
M[\lambda]$, $(v, a)\mapsto v\oor{\lambda} a$, satisfying:
\begin{align*}
(\partial v)\oor{\lambda} a&=-\lambda v\oor{\lambda} a,\qquad
v\oor{\lambda}(\partial a)=(\partial+\lambda)(v\oor{\lambda} a),\\
&(v\oor{\lambda} a)\oor{\lambda+\mu}b=v\oor{\lambda}(a\oo{\mu} b),
\end{align*}
for any $a, b\in A$ and $v\in M$. We denote it by $(M, \triangleleft)$.

A {\rm (conformal) $A$-bimodule} is a triple $(M, \triangleright, \triangleleft)$
such that $(M, \triangleright)$ is a left $A$-module, $(M, \triangleleft)$ is a
right $A$-module, and they satisfy
\begin{eqnarray*}
(a\ool{\lambda} v)\oor{\lambda+\mu}b=a\ool{\lambda}(v\oor{\mu}b),
\end{eqnarray*}
for any $a, b\in A$ and $v\in M$.
\end{defi}

Let $A$ be an associative conformal algebra. Define two bilinear maps
$\triangleright^{A},\; \triangleleft^{A}:\; A\otimes A\rightarrow A$ by
$a\oll{\lambda}{A} b=a\oo{\lambda} b$ and
$b\orr{\lambda}{A} a=b\oo{\lambda}a$ for all $a$, $b\in A$.
Then $(A, \triangleright^{A}, \triangleleft^{A})$ is a bimodule of $A$,
and it is called the regular bimodule of $A$.

Next, we recall the Hochschild type cohomology of an associative conformal algebra.
In 1999, Bakalov, Kac and Voronov first gave the definition of Hochschild
cohomology of associative conformal algebras \cite{BKV}. This definition is a
conformal analogue of Hochschild cohomology of associative algebras.
In \cite{DK}, for the case of Lie conformal algebras, the definition was improved by
taking $n-1$ variables. Following this idea, we define the Hochschild cohomology for
an associative conformal algebra $A$ by a bimodule $M$. We denote
$\mathcal{C}^{0}(A, M)=M/\partial M$, $\mathcal{C}^{1}(A, M)=\Hom_{\cb[\partial]}
(A, M)$, the set of $\cb[\partial]$-module homomorphisms from $A$ to $M$, and for
$n\geq2$, the space of $n$-cochains $\mathcal{C}^{n}(A, M)$ consists of all conformal
sesquilinear maps from $A^{\otimes n}$ to $M[\lambda_{1},\dots, \lambda_{n-1}]$,
i.e., the $\cb$-linear maps
$$
\varphi_{\lambda_{1},\dots, \lambda_{n-1}}:\;  A^{\otimes n}\longrightarrow
M[\lambda_{1},\dots, \lambda_{n-1}],
$$
such that
$$
\varphi_{\lambda_{1},\dots, \lambda_{n-1}}(a_{1},\dots, \partial a_{i},\dots, a_{n})
=-\lambda_{i}\varphi_{\lambda_{1},\dots, \lambda_{n-1}}(a_{1},\dots, a_{n}),
$$
for $i=1,2,\cdots,n-1$ and
$$
\varphi_{\lambda_{1},\dots, \lambda_{n-1}}(a_{1},\dots,  a_{n-1},\partial a_{n})
=(\partial+\lambda_{1}+\dots+\lambda_{n-1})\varphi_{\lambda_{1},\dots, \lambda_{n-1}}
(a_{1},\dots, a_{n}).
$$
The differentials are defined by $d_{0}: M/\partial M\rightarrow
\Hom_{\cb[\partial]}(A, M)$,
\begin{equation*}
d_{0}(v+\partial M)(a)=(a\ool{-\lambda-\partial}v-v\oor{\lambda}a)\mid_{\lambda=0},
\end{equation*}
and for $n\geq1$,
\begin{multline}\nonumber
\qquad d_{n}(\varphi)_{\lambda_{1},\dots, \lambda_{n}}(a_{1},\dots, a_{n+1})
=a_{1}\ool{\lambda_{1}}\varphi_{\lambda_{2},\dots, \lambda_{n}}
(a_{2},\dots, a_{n+1}) \\
+\sum_{i=1}^{n}(-1)^{i}\varphi_{\lambda_{1},\dots, \lambda_{i-1},
\lambda_{i}+\lambda_{i+1}, \lambda_{i+2}\dots,\lambda_{n}}
(a_{1},\dots, a_{i-1}, a_{i}\oo{\lambda_{i}} a_{i+1} , a_{i+2}\dots, a_{n+1}) \\
+(-1)^{n+1}\varphi_{\lambda_{1},\dots, \lambda_{n-1}}(a_{1},\dots,
a_{n})\oor{\lambda_{1}+\dots+\lambda_{n}} a_{n+1}.\qquad
\end{multline}
One can verify that the operator $d_{n}$ preserves the space of cochains and
$d_{n+1}\circ d_{n}=0$. The cochains of an associative conformal algebra $A$ with
coefficients in a bimodule $M$ form a complex $(\mathcal{C}^{\ast}(A, M),
d_{\ast})$, called the {\it Hochschild complex}. We denote the space of
$n$-cocycles by $\mathcal{Z}^{n}(A, M)=\{\varphi\in\mathcal{C}^{n}(A, M)\mid
d_{n}(\varphi)=0\}$, and the space of $n$-coboundaries by $\mathcal{B}^{n}
(A, M)=\{d_{n-1}(\varphi) \mid\varphi\in\mathcal{C}^{n-1}(A, M)\}$. The $n$-th
Hochschild cohomology of $A$ with coefficients in $M$ is defined by
\begin{equation*}
\HH^{n}(A, M)=\mathcal{Z}^{n}(A, M)/\mathcal{B}^{n}(A, M).
\end{equation*}
In particular, if $M=A$ as conformal bimodule, we denote $\HH^{n}(A):=\HH^{n}(A, A)$,
and call the $n$-th Hochschild cohomology of $A$. For the Hochschild cohomology of
an associative conformal algebra $A$, we have
\begin{itemize}
\item[(i)] $\Img d_{0}=\Inn(A, M)$, where $\Inn(A, M)=\{f_{v}\in\Hom_{\cb[\partial]}
           (A, M)\mid v\in M, f_{v}(a)=a\oo{-\partial}v-v\oo{0}a\}$.
\item[(ii)] $\Ker d_{1}=\Der(A, M)$, where $\Der(A, M)=\{f\in\Hom_{\cb[\partial]}
           (A, M)\mid f(a\oo{\lambda}b)=a\ool{\lambda}f(b)+f(a)\oor{\lambda}b\}$.
\item[(iii)] If $A$ as $\cb[\partial]$-module is projective, the equivalence classes
           of $\cb[\partial]$-split abelian extensions of $A$ by bimodule $M$
           is projective correspond bijectively to
           $\HH^{2}(A, M)$. For the details see \cite{Do,Lib}.
\end{itemize}

In \cite{HSZ}, for an associative conformal algebra $A$, we have shown that there is a
Gerstenhaber algebra structure on $\HH^{\bullet}(A)=\oplus_{n\geq0}\HH^{n}(A)$.
First, for the Hochschild complex $(\mathcal{C}^{\ast}(A, A), d_{\ast})$, we define the
Gerstenhaber bracket as following:
$$
[f, g]=f\bullet g-(-1)^{(m-1)(n-1)}g\bullet f
$$
for any $f\in\mathcal{C}^{m}(A, A)$ and $g\in\mathcal{C}^{n}(A, A)$, where
$f\bullet g=\sum_{i=0}^{m-1}(-1)^{(n-1)i}f\bullet_{i} g$, and
\begin{align*}
&(f\bullet_{i} g)_{\lambda_{0},\dots,\lambda_{m+n-3}}(a_{0}, a_{1},\dots, a_{m+n-2})\\
=&f_{\lambda_{0},\dots,\lambda_{i-1}, \lambda_{i}+\cdots+\lambda_{i+n},
\lambda_{i+n+1},\dots,\lambda_{m+n-3}}\big(a_{0},\dots, a_{i-1}, g_{\lambda_{i},
\dots,\lambda_{i+n-1}}(a_{i},\dots, a_{i+n}), a_{i+n+1}, \dots, a_{m+n-2}\big).
\end{align*}

\begin{pro}[\cite{HSZ}] \label{difLie}
Let $A$ be an associative conformal algebra. The triple $(\mathcal{C}^{\bullet+1}(A, A)
=\oplus_{n\geq0}\mathcal{C}^{n+1}(A, A),\; [-,-],\; d_{\ast})$ is a differential graded
Lie algebra.
\end{pro}

For a differential graded Lie algebra $(L=\oplus_{n\geq0}L^{i},\; [-,-],\; d)$,
recall that an element $\mathfrak{c}\in L^{1}$ is called a {\it Maurer Cartan element}
of this differential graded Lie algebra if $d(\mathfrak{c})+\frac{1}{2}[\mathfrak{c},
\mathfrak{c}]=0$. We denote the set of all the Maurer Cartan elements of
$(L=\oplus_{n\geq0}L^{i},\; [-,-],\; d)$ by ${\rm MC}(L)$.
If $L^{0}$ is abelian, there is a equivalence relation on ${\rm MC}(L)$ which is called
gauge equivalence relation (see \cite{Gol}). Let $\mathfrak{c}$ and $\tilde{\mathfrak{c}}$
be two elements in ${\rm MC}(L)$, they are called {\it gauge equivalent} if and only if there
exists $\xi\in L^{0}$ such that
$$
\tilde{\mathfrak{c}}=e^{\ad_{\xi}}(\mathfrak{c})-\frac{e^{\ad_{\xi}}-1}{\ad_{\xi}}d(\xi).
$$
The set of the gauge equivalence classes of ${\rm MC}(L)$ is denoted by $\mathcal{MC}(L)$.

Denote the element in $\mathcal{C}^{2}(A, A)$ corresponding to the conformal multiplication
by $\mathfrak{m}_{A}$, i.e., $\mathfrak{m}_{A}(a_{1}, a_{2})=a_{1}\oo{\lambda}a_{2}$
for $a_{1}, a_{2}\in A$. Then $\mathfrak{m}_{A}$ is a Maurer-Cartan element of
$(\mathcal{C}^{\bullet+1}(A, A),\; [-,-]$, $d_{\ast})$, and for any $f\in\mathcal{C}^{n}(A, A)$,
$d_{n}(f)=(-1)^{n-1}[\mathfrak{m}_{A}, f]$. Moreover, this bracket induces a degree $-1$
bracket $[-, -]$ on $\HH^{\bullet}(A)=\oplus_{n\geq0}\HH^{n}(A)$, such that there is a
Gerstenhaber algebra structure on $\HH^{\bullet}(A)$.  For the details see \cite{HSZ}.
Define $\bar{d}_{n}=(-1)^{n-1}d_{n}$ for $n\geq 0$, i.e., $\bar{d}_{n}(f)=[\mathfrak{m}_{A}, f]$
for any $f\in\mathcal{C}^{n}(A, A)$. Then we get a new differential graded Lie algebra
$(\mathcal{C}^{\bullet+1}(A, A),\; [-,-],\; \bar{d}_{\ast})$. In Section \ref{sec:Deli},
we will use the Maurer-Cartan elements of a subalgebra of this kind of differential
graded Lie algebra to characterize the non-abelian extensions of associative conformal algebras.

%%%%%%%%%%%%%%%%%%%%%%%%%%%%%%%%%%%%%%%%%%%%%%%%%%%%%%%%%%%%%%%%%%%%%%%%%%%%%%%%%
%    section  3   crossed module
%%%%%%%%%%%%%%%%%%%%%%%%%%%%%%%%%%%%%%%%%%%%%%%%%%%%%%%%%%%%%%%%%%%%%%%%%%%%%%%%%%%%%%

\section{Crossed modules and 2-term strongly homotopy associative conformal
algebras}\label{sec:cross}

In this section, we introduce the notions of crossed module of associative conformal algebras
and 2-term strongly homotopy associative conformal algebras. We will discuss the
relationship between them and the 3-th Hochschild cohomology of associative conformal
algebras with coefficient in bimodules.

\begin{defi}\label{def:2-t}
A {\rm 2-term strongly homotopy associative conformal algebra} consists of a complex of
$\cb[\partial]$-modules $A_{1}\xrightarrow{\fd}A_{0}$, conformal sesquilinear maps
$\mathfrak{m}^{2}: A_{i}\times A_{j}\rightarrow A_{i+j}[\lambda]$, $i, j\in\{0, 1\}$,
and a sesquilinear map $\mathfrak{m}^{3}: A_{0}\times A_{0}\times A_{0}
\rightarrow A_{1}[\lambda_{1}, \lambda_{2}]$, such that for any $a, a_{1}, a_{2}, a_{3},
a_{4}\in A_{0}$ and $b, b_{1}, b_{2}\in A_{1}$,
\begin{align}
\mathfrak{m}^{2}_{\lambda}(b_{1}, b_{2})&=0,\nonumber\\
\fd(\mathfrak{m}^{2}_{\lambda}(a, b))&=\mathfrak{m}^{2}_{\lambda}(a, \fd(b)),\label{2-t1}\\
\fd(\mathfrak{m}^{2}_{\lambda}(b, a))&=\mathfrak{m}^{2}_{\lambda}(\fd(b), a),\label{2-t2}\\
\mathfrak{m}^{2}_{\lambda}(\fd(b_{1}), b_{2})
&=\mathfrak{m}^{2}_{\lambda}(b_{1}, \fd(b_{2})),\label{2-t3}\\
\fd(\mathfrak{m}^{3}_{\lambda_{1}, \lambda_{2}}(a_{1}, a_{2}, a_{3}))
&=\mathfrak{m}^{2}_{\lambda_{1}+\lambda_{2}}
(\mathfrak{m}^{2}_{\lambda_{1}}(a_{1}, a_{2}), a_{3})
-\mathfrak{m}^{2}_{\lambda_{1}}(a_{1},
\mathfrak{m}^{2}_{\lambda_{2}}(a_{2}, a_{3})),\label{2-t4}\\
\mathfrak{m}^{3}_{\lambda_{1}, \lambda_{2}}(\fd(b_{1}), a_{2}, a_{3})
&=\mathfrak{m}^{2}_{\lambda_{1}+\lambda_{2}}
(\mathfrak{m}^{2}_{\lambda_{1}}(b_{1}, a_{2}), a_{3})
-\mathfrak{m}^{2}_{\lambda_{1}}(b_{1},
\mathfrak{m}^{2}_{\lambda_{2}}(a_{2}, a_{3})),\label{2-t5}\\
\mathfrak{m}^{3}_{\lambda_{1}, \lambda_{2}}(a_{1}, \fd(b_{2}), a_{3})
&=\mathfrak{m}^{2}_{\lambda_{1}+\lambda_{2}}(
\mathfrak{m}^{2}_{\lambda_{1}}(a_{1}, b_{2}), a_{3})
-\mathfrak{m}^{2}_{\lambda_{1}}(a_{1},
\mathfrak{m}^{2}_{\lambda_{2}}(b_{2}, a_{3})),\label{2-t6}\\
\mathfrak{m}^{3}_{\lambda_{1}, \lambda_{2}}(a_{1}, a_{2}, \fd(b_{3}))
&=\mathfrak{m}^{2}_{\lambda_{1}+\lambda_{2}}(
\mathfrak{m}^{2}_{\lambda_{1}}(a_{1}, a_{2}), b_{3})
-\mathfrak{m}^{2}_{\lambda_{1}}(a_{1},
\mathfrak{m}^{2}_{\lambda_{2}}(a_{2}, b_{3})),\label{2-t7}\\
\mathfrak{m}^{2}_{\lambda_{1}+\lambda_{2}+\lambda_{3}}
(\mathfrak{m}^{3}_{\lambda_{1}, \lambda_{2}}(a_{1}, a_{2}, a_{3}), &\, a_{4})
+\mathfrak{m}^{2}_{\lambda_{1}}(a_{1},
\mathfrak{m}^{3}_{\lambda_{2}, \lambda_{3}}(a_{2}, a_{3}, a_{4}))\label{2-t8}\\
=\mathfrak{m}^{3}_{\lambda_{1}+\lambda_{2}, \lambda_{3}}
(\mathfrak{m}^{2}_{\lambda_{1}}(a_{1}, a_{2}), a_{3}, a_{4})
&-\mathfrak{m}^{3}_{\lambda_{1}, \lambda_{2}+\lambda_{3}}(a_{1},
\mathfrak{m}^{2}_{\lambda_{2}}(a_{2}, a_{3}), a_{4})
+\mathfrak{m}^{3}_{\lambda_{1}, \lambda_{2}}(a_{1}, a_{2},
\mathfrak{m}^{2}_{\lambda_{3}}(a_{3}, a_{4})).\nonumber
\end{align}
We denote this 2-term strongly homotopy associative conformal algebra by
$(A_{1}\xrightarrow{\fd}A_{0}, \mathfrak{m}^{2}, \mathfrak{m}^{3})$, and it is said
to be {\rm skeletal} if $\fd=0$, is said to be {\rm strict} if $\mathfrak{m}^{3}=0$.
\end{defi}

The strongly homotopy associative algebras, in particular, the 2-term strongly homotopy
associative algebras were studied in \cite{Das, Mes}. The 2-term strongly homotopy
associative conformal algebras can be viewed as a conformal analogue of 2-term
strongly homotopy associative algebras. Following, we will simply write ``2-term strongly
homotopy associative conformal algebra" by ``2-term SHAC-algebra".

\begin{ex}
Let $A$ be an associative conformal algebra and $f: M_{1}\rightarrow M_{0}$
be a homomorphism between bimodules over $A$. We consider the 2-term chain complex
of $\cb[\partial]$-modules: $M_{1}\xrightarrow{\fd}A\oplus M_{0}$. It can be given
the structure of a strict 2-term SHAC-algebra by setting
\begin{align*}
\mathfrak{m}^{2}_{\lambda}((a, u), (b, v))&=(a\oo{\lambda}b,\;
a\ool{\lambda}v+u\oor{\lambda}b),\qquad
\mathfrak{m}^{2}_{\lambda}((a, u), v')=a\ool{\lambda}v',\\
\mathfrak{m}^{2}_{\lambda}(v', (a, u))&=v'\oor{\lambda}a,\qquad\qquad\qquad\qquad
\qquad\qquad\;\; \fd(v')=(0, f(v')),
\end{align*}
and $\mathfrak{m}^{3}=0$, for any $a, b\in A$, $u, v\in M_{0}$ and $v'\in M_{1}$.
\end{ex}

Let $A=(A_{1}\xrightarrow{\fd}A_{0}, \mathfrak{m}^{2}, \mathfrak{m}^{3})$ and
$A'=(A'_{1}\xrightarrow{\fd'}A'_{0}, \mathfrak{m}'^{2}, \mathfrak{m}'^{3})$
be two 2-term SHAC-algebras. A {\it morphism} $f=(f^{0}, f^{1}, f^{2}):
A\rightarrow A'$ consists of a chain map $(f^{1}, f^{0})$ from $A_{1}\xrightarrow{\fd}A_{0}$
to $A'_{1}\xrightarrow{\fd'}A'_{0}$, and a conformal sesquilinear map $f^{2}:
A_{0}\times A_{0}\rightarrow A'_{1}[\lambda]$ such that for any
$a, a_{1}, a_{2}, a_{3}\in A_{0}$ and $b\in A_{1}$,
\begin{align*}
\fd'(f^{2}_{\lambda}(a_{1}, a_{2}))&=f^{0}(\mathfrak{m}^{2}_{\lambda}(a_{1}, a_{2}))
-\mathfrak{m}'^{2}_{\lambda}(f^{0}(a_{1}), f^{0}(a_{2})),\\
f^{2}_{\lambda}(a, \fd(b))&=f^{1}(\mathfrak{m}^{2}_{\lambda}(a, b))
-\mathfrak{m}'^{2}_{\lambda}(f^{0}(a), f^{1}(b)),\\
f^{2}_{\lambda}(\fd(b), a)&=f^{1}(\mathfrak{m}^{2}_{\lambda}(b, a))
-\mathfrak{m}'^{2}_{\lambda}(f^{1}(b), f^{0}(a)),\\
f^{2}_{\lambda_{1}+\lambda_{2}}(\mathfrak{m}^{2}_{\lambda_{1}}(a_{1}, a_{2}), a_{3})
-f^{2}_{\lambda_{1}}&(a_{1}, \mathfrak{m}^{2}_{\lambda_{2}}(a_{2}, a_{3}))
+\mathfrak{m}'^{2}_{\lambda_{1}}(f^{0}(a_{1}), f^{2}_{\lambda_{2}}(a_{2},a_{3}))\\
=f^{1}(\mathfrak{m}^{3}_{\lambda_{1}, \lambda_{2}}(a_{1}, a_{2}, a_{3}))
-\mathfrak{m}'^{3}_{\lambda_{1}, \lambda_{2}}&(f^{0}(a_{1}), f^{0}(a_{2}), f^{0}(a_{3}))
+\mathfrak{m}'^{2}_{\lambda_{1}+\lambda_{2}}(f^{2}_{\lambda_{1}}(a_{1}, a_{2}), f^{0}(a_{3})).
\end{align*}
If $f=(f^{0}, f^{1}, f^{2}): A\rightarrow A'$  and $g=(g^{0}, g^{1}, g^{2}): A'\rightarrow A''$
be two morphism of 2-term SHAC-algebras, their composition $g\circ f: A\rightarrow A''$ is
defined by $(g\circ f)^{0}=g^{0}\circ f^{0}$, $(g\circ f)^{1}=g^{1}\circ f^{1}$, and
$$
(g\circ f)^{2}_{\lambda}(a_{1}, a_{2})=g^{2}_{\lambda}(f^{0}(a_{1}), f^{0}(a_{2}))
+g^{1}(f^{2}_{\lambda}(a_{1}, a_{2})),
$$
for any $a_{1}, a_{2}\in A_{0}$. For any 2-term SHAC-algebra $A$, there is an
identity morphism $\id_{A}=(\id_{A_{0}}, \id_{A_{1}}, 0)$.

\begin{pro}\label{pro:cat}
The collection of 2-term SHAC-algebras and morphisms between them form a category.
\end{pro}

Here we mainly consider skeletal 2-term SHAC-algebras, strict 2-term SHAC-algebras,
and the relationship between cohomology of associative conformal algebras and them.

\begin{defi}\label{def:2-t-e}
Let $A=(A_{1}\xrightarrow{0}A_{0}, \mathfrak{m}^{2}, \mathfrak{m}^{3})$ and
$A'=(A_{1}\xrightarrow{0}A_{0}, \mathfrak{m}'^{2}, \mathfrak{m}'^{3})$ be two
skeletal 2-term SHAC-algebras on the same chain complex of $\cb[\partial]$-modules.
They are said to be {\rm equivalent} if $\mathfrak{m}^{2}=\mathfrak{m}'^{2}$ and
there exists a conformal sesquilinear map $\sigma: A_{0}\times A_{0}\rightarrow A_{1}[\lambda]$
such that
$$
\mathfrak{m}'^{3}=\mathfrak{m}^{3}+d_{2}(\sigma),
$$
where $d_{2}$ is the differential in the Hochschild type cohomology complex.
\end{defi}

The following theorem characterizes the skeletal 2-term SHAC-algebras, and it also
gives a method to construct skeletal 2-term SHAC-algebras.

\begin{thm} \label{sk-alg}
There is a 1-1 correspondence between equivalence classes of skeletal
2-term SHAC-algebras and 3-th Hochschild cohomologies of associative conformal algebras
with coefficient in bimodules.
\end{thm}

\begin{proof}
Let $A=(A_{1}\xrightarrow{0}A_{0}, \mathfrak{m}^{2}, \mathfrak{m}^{3})$ be a
skeletal 2-term SHAC-algebra. Then $(A_{0}, \mathfrak{m}^{2})$ is an associative
conformal algebra, $\mathfrak{m}^{2}: A_{0}\times A_{1}\rightarrow A_{1}[\lambda]$
and $\mathfrak{m}^{2}: A_{1}\times A_{0}\rightarrow A_{1}[\lambda]$ give an $A_{0}$-bimodule
structure on $A_{1}$. By the equation (\ref{2-t8}) in Definition \ref{def:2-t},
we get $\mathfrak{m}^{3}$ is a 3-cocycle of $A_{0}$ with coefficients in bimodule $A_{1}$.

Conversely, given a quintuple $(A, M, \triangleright, \triangleleft, \zeta)$,
where $A$ is an associative conformal algebra, $(M, \triangleright, \triangleleft)$
is a bimodule over $A$, and $\zeta$ is a 3-cocycle of $A$ with coefficients in bimodule $M$.
We denote $A_{0}:=A$, $A_{1}:=M$, and define conformal sesquilinear maps
$\mathfrak{m}^{2}: A_{i}\times A_{j}\rightarrow A_{i+j}[\lambda]$ and $\mathfrak{m}^{3}:
A_{0}\times A_{0}\times A_{0}\rightarrow A_{1}[\lambda_{1}, \lambda_{2}]$ by
\begin{align*}
\mathfrak{m}^{2}_{\lambda}(a, b)&=a\oo{\lambda}b,\qquad\qquad\qquad\;\,
\mathfrak{m}^{2}_{\lambda}(a, v)=a\ool{\lambda}v,\\
\mathfrak{m}^{2}_{\lambda}(v, a)&=v\oor{\lambda}a,\qquad\qquad
\mathfrak{m}^{3}_{\lambda_{1}, \lambda_{2}}(a, b, c)=\zeta_{\lambda_{1}, \lambda_{2}}(a, b, c),
\end{align*}
for any $a, b, c\in A$ and $v\in M$. Then it is easy to verify that
$(A_{1}\xrightarrow{0}A_{0}, \mathfrak{m}^{2}, \mathfrak{m}^{3})$ is a skeletal
2-term SHAC-algebra, and the above two correspondences are inverses to each other.

If $A=(A_{1}\xrightarrow{0}A_{0}, \mathfrak{m}^{2}, \mathfrak{m}^{3})$ and
$A'=(A_{1}\xrightarrow{0}A_{0}, \mathfrak{m}'^{2}, \mathfrak{m}'^{3})$
are two equivalent skeletal 2-term SHAC-algebras, where the equivalence is given
by a conformal sesquilinear map $\sigma: A_{0}\times A_{0}\rightarrow A_{1}[\lambda]$,
i.e., $\mathfrak{m}'^{3}=\mathfrak{m}^{3}+d_{2}(\sigma)$. Then $[\mathfrak{m}^{3}]=
[\mathfrak{m}'^{3}]$, where $[\mathfrak{m}^{3}]$ means the equivalence class of
$\mathfrak{m}^{3}$ in $\HH^{3}(A_{0}, A_{1})$, where the $A_{0}$-bimodule structure on
$A_{1}$ is given by $\mathfrak{m}^{2}$. We denote by $\mathcal{S}kel(A_{0}, A_{1},
\mathfrak{m}^{2})$ the set of equivalent classes of skeletal 2-term SHAC-algebra
structures on $A_{1}\xrightarrow{0}A_{0}$ with same $\mathfrak{m}^{2}$. Then there is a map
$$
\Upsilon:\quad \mathcal{S}kel(A_{0}, A_{1}, \mathfrak{m}^{2})\longrightarrow
\HH^{3}(A_{0}, A_{1}),
\qquad [(A_{1}\xrightarrow{0}A_{0}, \mathfrak{m}^{2}, \mathfrak{m}^{3})]\mapsto
[\mathfrak{m}^{3}].
$$
Conversely, let $A_{0}$ be an associative conformal algebra and $A_{1}$ be an $A_{0}$-bimodule.
It is easy to see that any two representatives of an element
in $\HH^{3}(A_{0}, A_{1})$ give two equivalent skeletal 2-term SHAC-algebras
on $A_{1}\xrightarrow{0}A_{0}$ with same $\mathfrak{m}^{2}$, where $\mathfrak{m}^{2}$
is given by the module action of $A_{0}$ on $A_{1}$ and the multiplication of $A_{0}$.
Thus, the map $\Upsilon$ induces a 1-1 correspondence between equivalence classes
of skeletal 2-term SHAC-algebras and the elements in 3-th Hochschild cohomology of associative
conformal algebras with coefficient in bimodules.
\end{proof}

Next, we consider the crossed module of associative conformal algebras.

\begin{defi}\label{def:act}
Let $(X,\; \cdot\ooo{\lambda}{X}\cdot)$, $(Y,\; \cdot\ooo{\lambda}{Y}\cdot)$ be two
associative conformal algebras. If there are two linear maps
$\triangleright: X\times Y\rightarrow Y[\lambda]$ and $\triangleleft: Y\times X\rightarrow
Y[\lambda]$ such that $(Y, \triangleright, \triangleleft)$ is a bimodule over $X$ and
for any $x\in X$ and $y_{1}, y_{2}\in Y$,
\begin{align*}
(x\ool{\lambda} y_{1})\ooo{\lambda+\mu}{Y}y_{2}
&=x\ool{\lambda}(y_{1}\ooo{\mu}{Y}y_{2}),\\
(y_{1}\oor{\lambda}x)\ooo{\lambda+\mu}{Y}y_{2}
&=y_{1}\ooo{\lambda}{Y}(x\ool{\mu}y_{2}),\\
(y_{1}\ooo{\lambda}{Y}y_{2})\oor{\lambda+\mu}x
&=y_{1}\ooo{\lambda}{Y}(y_{2}\oor{\mu}x),
\end{align*}
then we call {\rm $X$ acts on $Y$} by $(\triangleright, \triangleleft)$,
or simply $X$ acts on $Y$.
\end{defi}

Following, we will simply denote $x_{1}\ooo{\lambda}{X}x_{2}$ by $x_{1}\oo{\lambda}x_{2}$
and denote $y_{1}\ooo{\lambda}{Y}y_{2}$ by $y_{1}\oo{\lambda}y_{2}$ for any
$x_{1}, x_{2}\in X$ and $y_{1}, y_{2}\in Y$, if it does not affect understanding.

\begin{defi}\label{def:cross}
Let $X$, $Y$ be two associative conformal algebras. A {\rm crossed module of associative
conformal algebras} is a quintuple $(X, Y, \rho, \triangleright, \triangleleft)$, where
$(\triangleright, \triangleleft)$ is an action of $X$ on $Y$, $\rho: Y\rightarrow X$
is a homomorphism of $\cb[\partial]$-modules such that for any $x\in X$ and
$y, y_{1}, y_{2}\in Y$,
\begin{align*}
\rho(x\ool{\lambda} y)&=x\oo{\lambda}\rho(y),\\
\rho(y\oor{\lambda} x)&=\rho(y)\oo{\lambda}x,\\
\rho(y_{1})\ool{\lambda} y_{2}=y_{1}&\oo{\lambda}y_{2}=y_{1}\oor{\lambda}\rho(y_{2}).
\end{align*}
\end{defi}

The crossed module of groups, associative algebras, Lie algebras and Hopf algebras
have been extensively studied.
See literature \cite{Wag} for details. The crossed module of associative conformal
algebras can be viewed as a conformal analogue of crossed module of associative algebras.
We will just say ``crossed module" instead of ``crossed module of associative conformal
algebras" except when emphasis is needed, and we remark that from the above crossed
module conditions one get $\rho$ is an associative conformal algebra homomorphism:
$\rho(y_{1}\oo{\lambda} y_{2})=\rho(\rho(y_{1})\ool{\lambda} y_{2})
=\rho(y_{1})\oo{\lambda}\rho(y_{2})$.

\begin{ex}
Let $I$ be a two sided ideal of an associative conformal algebra $A$. Then
$(I, A, \iota, \triangleright, \triangleleft)$ is a crossed module, where $\iota$ is the
canonical inclusion map and the action $(\triangleright, \triangleleft)$ is given by the
multiplication of $A$. In particular, if $I=0$ or $I=A$, then $(0, A, 0, 0, 0)$ and
$(A, A, \id, \cdot\oo{\lambda}\cdot, \cdot\oo{\lambda}\cdot)$ are crossed modules.
\end{ex}

\begin{thm} \label{cross-mod}
There is a 1-1 correspondence between strict 2-term SHAC-algebras and
crossed modules of associative conformal algebras.
\end{thm}

\begin{proof}
Let $A=(A_{1}\xrightarrow{\fd}A_{0}, \mathfrak{m}^{2}, \mathfrak{m}^{3})$ be a
strict 2-term SHAC-algebra. Let $X:=A_{0}$. Then, by the equation (\ref{2-t4}) in
Definition \ref{def:2-t}, $X$ with the multiplication
$x_{1}\oo{\lambda}x_{2}=\mathfrak{m}^{2}_{\lambda}(x_{1}, x_{2})$, $x_{1}, x_{2}\in X$,
is an associative conformal algebra. Set $Y:=A_{1}$ and
$y_{1}\oo{\lambda}y_{2}=\mathfrak{m}^{2}_{\lambda}(\fd(y_{1}), y_{2})
=\mathfrak{m}^{2}_{\lambda}(y_{1}, \fd(y_{2}))$, Then for any $y_{1}, y_{2}, y_{3}\in Y$,
since $\mathfrak{m}^{3}=0$, by the equation (\ref{2-t7}), we get
\begin{align*}
&\; y_{1}\oo{\lambda}(y_{2}\oo{\mu} y_{3})-(y_{1}\oo{\lambda} y_{2})\oo{\lambda+\mu} y_{3}\\
=&\; \mathfrak{m}^{2}_{\lambda}(\fd(y_{1}), \mathfrak{m}^{2}_{\mu}(\fd(y_{2}), y_{3}))
-\mathfrak{m}^{2}_{\lambda+\mu}(\mathfrak{m}^{2}_{\lambda}(\fd(y_{1}), \fd(y_{2})), y_{3})\\
=&\; \mathfrak{m}^{3}_{\lambda, \mu}(\fd(y_{1}), \fd(y_{2}), \fd(y_{3}))\\
=&\;0.
\end{align*}
This means that $Y$ is an associative conformal algebra. Moreover, by equations
(\ref{2-t5})-(\ref{2-t7}), we obtain that $\triangleright:=\mathfrak{m}^{2}:
X\times Y\rightarrow Y[\lambda]$ and $\triangleleft:=\mathfrak{m}^{2}:
Y\times X\rightarrow Y[\lambda]$ give an action of $X$ on $Y$.
Finally, take $\rho=\fd: Y\rightarrow X$.
Then $\rho$ is a homomorphism of $\cb[\partial]$-modules and the equations
(\ref{2-t1})-(\ref{2-t3}) are just the conditions in Definition \ref{def:cross}.
Thus the quintuple $(X, Y, \rho, \triangleright, \triangleleft)$ is a crossed module.

Conversely, given a crossed module $(X, Y, \rho, \triangleright, \triangleleft)$,
we can construct a strict 2-term SHAC-algebra as follows. Set $A_{0}=X$, $A_{1}=Y$,
$\fd=\rho$, and conformal sesquilinear maps $\mathfrak{m}^{2}: A_{i}\times A_{j}
\rightarrow A_{i+j}[\lambda]$ by
\begin{align*}
\mathfrak{m}^{2}_{\lambda}(x_{1}, x_{2})&=x_{1}\oo{\lambda}x_{2},\qquad\qquad\qquad\;
\mathfrak{m}^{2}_{\lambda}(x, y)=x\ool{\lambda}y,\\
\mathfrak{m}^{2}_{\lambda}(y, x)&=y\oor{\lambda}x,\qquad\qquad\qquad
\mathfrak{m}^{2}_{\lambda}(y_{1}, y_{2})=0,
\end{align*}
for any $x, x_{1}, x_{2}\in X$ and $y, y_{1}, y_{2}\in Y$. Then it is easy to verify
that $(A_{1}\xrightarrow{\fd}A_{0}, \mathfrak{m}^{2}, \mathfrak{m}^{3})$ is a
strict 2-term SHAC-algebra. Thus we get a 1-1 correspondence between strict 2-term
SHAC-algebras and crossed modules of associative conformal algebras.
\end{proof}

At the end of this section, we relate the 3-th Hochschild cohomology of associative conformal
algebras by means of a particular kind of crossed modules.

\begin{defi}\label{def:cross-ext}
Let $A$ be an associative conformal algebra and $(M, \blacktriangleright, \blacktriangleleft)$
be a bimodule over $A$. A {\rm crossed extension} of $A$ by $M$ is an exact sequence
of associative conformal algebras
$$
\mathcal{S}:\qquad \xymatrix{0\ar[r]& M\ar[r]^{\alpha}&Y\ar[r]^{\beta}
&X\ar[r]^{\gamma}&A\ar[r]&0,}
$$
such that there is an action $(\triangleright, \triangleleft)$ of $X$ on $Y$ and
$(X, Y, \beta, \triangleright, \triangleleft)$ is a crossed module, where the multiplication
of associative conformal algebra $M$ is trivial.

The crossed extension $\mathcal{S}$ is said to be {\rm split} if there are
$\cb[\partial]$-module homomorphisms $\varrho: A\rightarrow X$ and $\varsigma:
\Img(g)\rightarrow Y$ such that $\gamma\circ\varrho=\id_{A}$,
$\beta\circ\varsigma=\id_{\Img(\beta)}$. Here $\varrho$ is called a section of $\mathcal{S}$.
\end{defi}

Given a crossed module $(X, Y, \rho, \triangleright, \triangleleft)$, since $\rho$
is a homomorphism of associative conformal algebras, there is
an exact sequence of associative conformal algebras
$$
\mathcal{S}':\qquad \xymatrix{0\ar[r]& M\ar[r]^{\iota}&Y\ar[r]^{\rho}&X\ar[r]^{\pi}&A\ar[r]&0,}
$$
where $M=\Ker(\rho)$ and $A=\Cok(\rho)$. Note that for any
$u, v\in M$, $u\oo{\lambda}v=\rho(u)\ool{\lambda}v=0$ since $\rho(u)=0$, we get the
multiplication of $M$ is trivial. If there is a $\cb[\partial]$-module
homomorphism $\varrho$ such that $\pi\circ\varrho=\id_{A}$, there is an
$A$-bimodule structure on $M$ as follows. For any $a\in A$ and $v\in M$,
we define $a\blacktriangleright_{\lambda}v:=\varrho(a)\ool{\lambda}m$
and $v\blacktriangleleft_{\lambda}a:=v\oor{\lambda}\varrho(a)$, then one can check that
$(M, \blacktriangleright, \blacktriangleleft)$ is an $A$-bimodule which is called
the $A$-bimodule structure on $M$ induced by $\varrho$, and it
does not depend on the choice of $\varrho$.

Following, each crossed extension of $A$ by $M$ refers to a crossed extension such that
the $A$-module structure on $M$ coincides with the $A$-bimodule structure
induced by the section.

\begin{defi}\label{def:equ-cross-ext}
Let $\xymatrix@C=0.4cm{0\ar[r]& M\ar[r]^{\alpha}&Y\ar[r]^{\beta}&X\ar[r]^{\gamma}
&A\ar[r]&0}$ and $\xymatrix@C=0.4cm{0\ar[r]& M\ar[r]^{\alpha'}&Y'\ar[r]^{\beta'}
&X'\ar[r]^{\gamma'}&A\ar[r]&0}$ be two split crossed extensions of $A$ by $M$.
If there are homomorphisms of associative conformal algebras $\varphi: Y\rightarrow Y'$
and $\psi: X\rightarrow X'$ such that the following diagram is commutative
$$
\xymatrix{0\ar[r] & M\ar[r]^{\alpha}\ar@{=}[d] & Y\ar[r]^{\beta}\ar[d]^{\varphi} &
X\ar[r]^{\gamma}\ar[d]^{\psi} & A\ar[r]\ar@{=}[d]&0\;\\
0\ar[r] & M\ar[r]^{\alpha'} & Y'\ar[r]^{\beta'} & X'\ar[r]^{\gamma'} & A\ar[r]&0,}
$$
we call $(\varphi, \psi)$ is a {\rm morphism of crossed extensions} of $A$ by $M$.

This two crossed extensions of $A$ by $M$ are said to be {\rm equivalent} if there is a
morphism from one to the other. Denote by $\Ext^{2}(A, M)$ the set of
equivalence classes of split crossed extensions of $A$ by $M$.
\end{defi}

\begin{pro}\label{pro: cross-coh}
For any associative conformal algebra $A$ and a bimodule $M$ over $A$, there is a canonical
map
$$
\Theta:\quad \Ext^{2}(A, M)\longrightarrow\HH^{3}(A, M).
$$
\end{pro}

\begin{proof}
Given a split crossed extension $\mathcal{S}:\xymatrix@C=0.4cm{0\ar[r]& M\ar[r]^{\alpha}
&Y\ar[r]^{\beta}&X\ar[r]^{\gamma}&A\ar[r]&0}$ with a section $\varrho: A\rightarrow X$ and
$\cb[\partial]$-module homomorphism $\varsigma: \Img(\beta)\rightarrow Y$ such that
$\beta\circ\varsigma=\id_{\Img(\beta)}$. We have $\varrho(a)\oo{\lambda}\varrho(b)
-\varrho(a\oo{\lambda}b)\in\Ker(\gamma)[\lambda]=\Img(\beta)[\lambda]$, for any $a, b\in A$.
Denote $g_{\lambda}(a, b)=\varsigma(\varrho(a)\oo{\lambda}\varrho(b)-\varrho(a\oo{\lambda}b))$
and define
\begin{align*}
f_{\lambda_{1}, \lambda_{2}}(a, b, c)
&=\varrho(a)\ool{\lambda_{1}}g_{\lambda_{2}}(b, c)
-g_{\lambda_{1}+\lambda_{2}}(a\oo{\lambda_{1}}b,\; c)\\
&\quad+g_{\lambda_{1}}(a,\; b\oo{\lambda_{2}}c)
-g_{\lambda_{1}}(a, b)\oor{\lambda_{1}+\lambda_{2}}\varrho(c),
\end{align*}
for any $a, b, c\in A$. Then one can check that $f$ is a conformal sesquilinear map, and
equation $\beta(f_{\lambda_{1}, \lambda_{2}, \lambda_{3}}(a, b, c))=0$ holds,
since $(X, Y, \beta, \triangleright, \triangleleft)$ is a crossed module.
Moreover, one can obtain $d_{3}(f)=0$ by routine calculations, where $d_{3}$
is a differential in the Hochschild cohomology complex of $A$ with coefficients in $M$.
We define $\Theta(\mathcal{S})=[f]$. Following we need to
show that the map $\Theta$ is well-defined.

First, we are going to show that $\Theta$ does not depend on the choice of $\varrho$.
Let $\bar{\varrho}$ be another section of $\mathcal{S}$ and
$\bar{f}$ be the 3-cocycle defined using $\bar{\varrho}$ instead of $\varrho$.
Then, there exists a $\cb[\partial]$-module homomorphism $\eta: A\rightarrow Y$ such that
$\beta\circ\eta=\bar{\varrho}-\varrho$. Denote by $\bar{g}_{\lambda}(a, b)=
\varsigma(\bar{\varrho}(a)\oo{\lambda}\bar{\varrho}(b)-\bar{\varrho}(a\oo{\lambda}b))$.
Then, by crossed modules properties, we have
\begin{align*}
(\bar{f}-f)_{\lambda_{1}, \lambda_{2}}(a, b, c)
&=\varrho(a)\ool{\lambda_{1}}(\bar{g}-g)_{\lambda_{2}}(b, c)
-(\bar{g}-g)_{\lambda_{1}+\lambda_{2}}(a\oo{\lambda_{1}}b,\; c)\\
&\quad+(\bar{g}-g)_{\lambda_{1}}(a,\; b\oo{\lambda_{2}}c)
-(\bar{g}-g)_{\lambda_{1}}(a, b)\oor{\lambda_{1}+\lambda_{2}}\varrho(c)\\
&\quad+\eta(a)\oo{\lambda_{1}}\bar{g}_{\lambda_{2}}(b, c)
-\bar{g}_{\lambda_{1}}(a, b)\oo{\lambda_{1}+\lambda_{2}}\eta(c),
\end{align*}
for any $a, b, c\in A$. Define conformal sesquilinear map $\tilde{g}:
A\times A\rightarrow Y[\lambda]$ by
$$
\tilde{g}_{\lambda}(a, b)=\bar{\varrho}(a)\ool{\lambda}\eta(b)
+\eta(a)\oor{\lambda}\bar{\varrho}(b)-\eta(a\oo{\lambda}b)-\eta(a)\oo{\lambda}\eta(b),
$$
for $a, b\in A$. Then it is easy to see that $\beta\circ\tilde{g}=\beta\circ(\bar{g}-g)$.
That is $\bar{g}-g-\tilde{g}\in\mathcal{C}^{2}(A, M)$. We denote
\begin{align*}
\tilde{f}_{\lambda_{1}, \lambda_{2}}(a, b, c)
&:=\varrho(a)\ool{\lambda_{1}}\tilde{g}_{\lambda_{2}}(b, c)
-\tilde{g}_{\lambda_{1}+\lambda_{2}}(a\oo{\lambda_{1}}b,\; c)\\
&\quad+\tilde{g}_{\lambda_{1}}(a,\; b\oo{\lambda_{2}}c)
-\tilde{g}_{\lambda_{1}}(a, b)\oor{\lambda_{1}+\lambda_{2}}\varrho(c)\\
&\quad+\eta(a)\oo{\lambda_{1}}\bar{g}_{\lambda_{2}}(b, c)
-\bar{g}_{\lambda_{1}}(a, b)\oo{\lambda_{1}+\lambda_{2}}\eta(c).
\end{align*}
Then, by the crossed modules properties, we obtain $\tilde{f}_{\lambda_{1},
\lambda_{2}}(a, b, c)=0$. Note that $\bar{f}-f-\tilde{f}=d_{2}(\bar{g}-g-\tilde{g})$,
we get $[\bar{f}-f]=[\tilde{f}]=0$ in $\HH^{3}(A, M)$, i.e., $[\bar{f}]=[f]$ in $\HH^{3}(A, M)$.

Let $\mathcal{S}':\xymatrix@C=0.4cm{0\ar[r]& M\ar[r]^{\alpha'}
&Y'\ar[r]^{\beta'}&X'\ar[r]^{\gamma'}&A\ar[r]&0}$ be another split crossed extension
of $A$ by $M$, and $(\varphi, \psi)$ be a morphism from $\mathcal{S}$ to $\mathcal{S}'$.
Denote by $(X', Y', \beta', \triangleright', \triangleleft')$ the corresponding crossed
module of $\mathcal{S}'$. Then $\varrho'=\psi\circ\varrho$ is a section of $\mathcal{S}'$,
and for any $v\in M$, $a\in A$, $\varrho'(a)\triangleright'v=\varrho(a)\triangleright v$,
$v\triangleleft\varrho(a)=v\triangleleft'\varrho'(a)$.
Let $\varsigma': \Img(\beta')\rightarrow Y'$ be an arbitrary $\cb[\partial]$-module
homomorphism such that $\beta'\circ\varsigma'=\id_{\Img(\beta')}$. As discussed above,
we can use $\varrho$ and $\varsigma$ to define $f$, use $\varrho'$ and $\varsigma'$ to
define $f'$. Now we define a conformal sesquilinear map $h: A\times A\rightarrow M[\lambda]$ by
$$
h_{\lambda}(a, b)=(\varphi\circ\varsigma-\varsigma'\circ\psi)
\Big(\varrho(a)\oo{\lambda}\varrho(b)-\varrho(a\oo{\lambda}b)\Big),
$$
for $a, b\in A$. Then, for any $a, b, c\in A$, we have
\begin{align*}
&\; (f-f')_{\lambda_{1}, \lambda_{2}}(a, b, c)\\
=&\; \varrho(a)\ool{\lambda_{1}}\varsigma\Big(\varrho(b)\oo{\lambda_{2}}\varrho(c)
-\varrho(b\oo{\lambda_{2}}c)\Big)
-\varsigma\Big(\varrho(a\oo{\lambda_{1}}b)\oo{\lambda_{1}+\lambda_{2}}\varrho(c)\Big)\\
&\;+\varsigma\Big(\varrho(a)\oo{\lambda_{1}}\varrho(b\oo{\lambda_{2}}c)\Big)
-\varsigma\Big(\varrho(a)\oo{\lambda_{1}}\varrho(b)
-\varrho(a\oo{\lambda_{1}}b)\Big)\oor{\lambda_{1}+\lambda_{2}}\varrho(c)\\
&\;-\varrho'(a)\oll{\lambda_{1}}{'}\varsigma'\Big(\varrho'(b)\oo{\lambda_{2}}\varrho'(c)
-\varrho'(b\oo{\lambda_{2}}c)\Big)
+\varsigma'\Big(\varrho'(a\oo{\lambda_{1}}b)\oo{\lambda_{1}+\lambda_{2}}\varrho'(c)\Big)\\
&\;-\varsigma'\Big(\varrho'(a)\oo{\lambda_{1}}\varrho'(b\oo{\lambda_{2}}c)\Big)
+\varsigma'\Big(\varrho'(a)\oo{\lambda_{1}}\varrho'(b)
-\varrho'(a\oo{\lambda_{1}}b)\Big)\orr{\lambda_{1}+\lambda_{2}}{'}\varrho'(c)\\
=&\; d_{2}(h)_{\lambda_{1}, \lambda_{2}}(a, b, c).
\end{align*}
That is to say, $[f]=[f']$ in $\HH^{3}(A, M)$. By the arbitrariness of $\varsigma'$,
we also get that $\Theta$ does not depend on the choice of $\varsigma$.
Thus, the map $\Theta$ is well-defined.
\end{proof}

As in the classical case, we would like to construct an isomorphism between
$\Ext^{2}(A, M)$ and $\HH^{3}(A, M)$. But now we cannot construct a
canonical example of split crossed extension for a given cohomology class
in $\HH^{3}(A, M)$. We need to further consider the conditions for map $\Theta$ to be
injective and surjective.

%%%%%%%%%%%%%%%%%%%%%%%%%%%%%%%%%%%%%%%%%%%%%%%%%%%%%%%%%%%%%%%%%%%%%%%%%%%%%%%%%
%    section  4   Non-abelian extension and cohomology
%%%%%%%%%%%%%%%%%%%%%%%%%%%%%%%%%%%%%%%%%%%%%%%%%%%%%%%%%%%%%%%%%%%%%%%%%%%%%%%%%%%%%%

\section{Non-abelian extension and cohomology}\label{sec:Non-abel}

In this section, we consider the non-abelian extensions of associative conformal algebras,
define a non-abelian cohomology group, and show that the non-abelian extensions can be
classified by the non-abelian cohomology group. First, we give the definition of
non-abelian extensions of associative conformal algebras.

\begin{defi}\label{def:nonabel}
Let $A$ and $B$ be two associative conformal algebras. A {\rm non-abelian extension}
$\mathcal{E}$ of $B$ by $A$ is a short exact sequence
$$
\mathcal{E}:\qquad \xymatrix{0\ar[r]& A\ar[r]^{\alpha}&E\ar[r]^{\beta}&B\ar[r]&0,}
$$
where $E$ is an associative conformal algebra, $\alpha, \beta$ are homomorphisms of
associative conformal algebras, and this sequence is split in the category of
$\cb[\partial]$-modules. Let $\mathcal{E}$ and $\mathcal{E}'$ be two extensions of $B$ by $A$.
They are called {\rm equivalent} if there exists a homomorphism of associative conformal
algebras $\theta: E\rightarrow E'$ such that the following diagram commutes
$$
\xymatrix{\mathcal{E}:\qquad 0\ar[r]& A\ar[r]^{\alpha}\ar@{=}[d]&E\ar[r]^{\beta}
\ar[d]^{\theta} &B\ar[r]\ar@{=}[d]&0\;\\ \mathcal{E}':\qquad 0\ar[r]&
A\ar[r]^{\alpha'}&E'\ar[r]^{\beta'}&B\ar[r]&0.}
$$
\end{defi}

\begin{rmk}\label{rmk1}
The non-abelian extension of associative conformal algebra defined here is also a
$\cb[\partial]$-split extension. In this regard, we can require associative conformal algebra
$B$ to be projective as $\cb[\partial]$-module. It is a generalization of $\cb[\partial]$-split
abelian extension of associative conformal algebras.
\end{rmk}

We now define a non-abelian cohomology and show that the non-abelian extensions are
classified by the non-abelian cohomology.

\begin{defi}\label{def:coho}
Let $A$ and $B$ be two associative conformal algebras. A {\rm non-abelian 2-cocycle} on $B$
with values in $A$ is a triplet $(\blacktriangleright, \blacktriangleleft, \chi)$ of
$\cb[\partial]$-bilinear maps $\chi: B\times B\rightarrow A[\lambda]$, $(b_{1}, b_{2})
\mapsto\chi_{\lambda}(b_{1}, b_{2})$, $\blacktriangleright:
B\times A\rightarrow A[\lambda]$, $(b, a)\mapsto b\obl{\lambda} a$, and $\blacktriangleleft:
A\times B\rightarrow A[\lambda]$, $(a, b)\mapsto a\obr{\lambda} b$, satisfying
the following properties:
\begin{align}
b_{1}\obl{\lambda}(b_{2}\obl{\mu}a)&=(b_{1}\oo{\lambda}b_{2})\obl{\lambda+\mu} a+
\chi_{\lambda}(b_{1}, b_{2})\oo{\lambda+\mu}a,  \label{coh1}\\
(a\obr{\lambda}b_{1})\obr{\lambda+\mu}b_{2}&=a\obr{\lambda}(b_{1}\oo{\mu}b_{2})+
a\oo{\lambda}\chi_{\mu}(b_{1}, b_{2}),  \label{coh2}\\
b_{1}\obl{\lambda}(a\obr{\mu}b_{2})&=(b_{1}\obl{\lambda}a)
\obr{\lambda+\mu}b_{2}, \label{coh3}\\
(a_{1}\oo{\lambda}a_{2})\obr{\lambda+\mu}b&=a_{1}\oo{\lambda}(a_{2}\obr{\mu} b),\label{coh4}\\
b\obl{\lambda}(a_{1}\oo{\mu}a_{2})&=(b\obl{\lambda}a_{1})\oo{\lambda+\mu}a_{2},\label{coh4'}\\
(a_{1}\obr{\lambda} b)\oo{\lambda+\mu}a_{2}&=a_{1}\oo{\lambda}(b\obl{\mu} a_{2}),
\label{coh4''}\\
b_{1}\obl{\lambda}\chi_{\mu}(b_{2}, b_{3})+\chi_{\lambda+\mu}(b_{1},\; b_{2}\oo{\mu}b_{3})
&=\chi_{\lambda+\mu}(b_{1}\oo{\lambda}b_{2},\; b_{3})
+\chi_{\lambda}(b_{1}, b_{2})\obr{\lambda+\mu} b_{3},\label{coh5}
\end{align}
for any $a, a_{1}, a_{2}\in A$ and $b, b_{1}, b_{2}, b_{3}\in B$.
One denotes by $\mathcal{Z}^{2}_{nab}(B, A)$ the set of theses cocycles.

Moreover, $(\blacktriangleright, \blacktriangleleft, \chi)$ and $(\bar{\blacktriangleright},
\bar{\blacktriangleleft}, \bar{\chi})$ are said to be {\rm equivalent} if there exists
a $\cb[\partial]$-module homomorphism $\delta: B \rightarrow A$ satisfying:
\begin{align}
b\obll{\lambda}a-b\obl{\lambda}a &=\delta(b)\oo{\lambda}a,  \label{coh6}\\
a\obrr{\lambda}b-a\obr{\lambda}b &=a\oo{\lambda}\delta(b),  \label{coh7}\\
\bar{\chi}_{\lambda}(b_{1}, b_{2})-\chi_{\lambda}(b_{1}, b_{2})
&=b_{1}\obll{\lambda}\delta(b_{2})-\delta(b_{1}\oo{\lambda}b_{2})
+\delta(b_{1})\obrr{\lambda}b_{2}-\delta(b_{1})\oo{\lambda}\delta(b_{2}),  \label{coh8}
\end{align}
for any $a\in A$ and $b, b_{1}, b_{2}\in B$. In this case, we denote
$(\blacktriangleright, \blacktriangleleft, \chi)\approx(\bar{\blacktriangleright},
\bar{\blacktriangleleft}, \bar{\chi})$. The {\rm non-abelian cohomology}
$\HH^{2}_{nab}(B, A)$ is the quotient of $\mathcal{Z}^{2}_{nab}(B, A)$ by this
equivalence relation.
\end{defi}

Let $A$ and $B$ be two associative conformal algebras.
Let $\mathcal{E}$ be a non-abelian extension of $B$ by $A$, i.e., there is a short
exact sequence
$$
\xymatrix{0\ar[r]& A\ar[r]^{\alpha}&E\ar[r]^{\beta}&B\ar[r]&0},
$$
which is split in the category of $\cb[\partial]$-modules. Thus, there is a
$\cb[\partial]$-module homomorphism $\gamma: B\rightarrow E$
such that $\beta\circ\gamma=\id_{B}$, which is called a section of $\mathcal{E}$.
We define $\cb[\partial]$-bilinear maps $\chi^{\gamma}: B\times B\rightarrow
A[\lambda]$, $\blacktriangleright^{\gamma}: B\times A\rightarrow A[\lambda]$ and
$\blacktriangleleft^{\gamma}: A\times B\rightarrow A[\lambda]$ by the multiplication
of $E$ as following:
\begin{align*}
b\bll{\lambda}{\gamma}a=\gamma(b)\oo{\lambda}a&, \qquad\qquad\qquad\quad
a\brr{\lambda}{\gamma}b=a\oo{\lambda}\gamma(b),\\
\chi^{\gamma}_{\lambda}(b_{1}, b_{2})&=\gamma(b_{1})\oo{\lambda}\gamma(b_{2})
-\gamma(b_{1}\oo{\lambda} b_{2}),
\end{align*}
for any $a\in A$ and $b, b_{1}, b_{2}\in B$.

\begin{lem} \label{lem: 2-cyc}
With the above notations, the triplet $(\blacktriangleright^{\gamma},
\blacktriangleleft^{\gamma}, \chi^{\gamma})$ is a non-abelian 2-cocycle on $B$ with
values in $A$.
\end{lem}

\begin{proof}
By direct calculation, for any $a\in A$ and $b_{1}, b_{2}\in B$, we have
\begin{align*}
b_{1}\bll{\lambda}{\gamma}(b_{2}\bll{\mu}{\gamma}a)
&=\gamma(b_{1})\oo{\lambda}(\gamma(b_{2})\oo{\mu}a)\\
&=\Big(\gamma(b_{1})\oo{\lambda}\gamma(b_{2})-\gamma(b_{1}\oo{\lambda}b_{2})\Big)
\oo{\lambda+\mu}a+\gamma(b_{1}\oo{\lambda}b_{2})\oo{\lambda+\mu}a\\
&=(b_{1}\oo{\lambda}b_{2})\bll{\lambda+\mu}{\gamma} a+
\chi^{\gamma}_{\lambda}(b_{1}, b_{2})\oo{\lambda+\mu}a.
\end{align*}
This means the equation (\ref{coh1}) in Definition \ref{def:coho} is satisfied.
Similarly, one can check that the equations (\ref{coh2}-\ref{coh4''}) are satisfied too.
Considering the equation (\ref{coh5}), we have, for any $b_{1}, b_{2}, b_{3}\in B$,
\begin{align*}
&b_{1}\bll{\lambda}{\gamma}\chi^{\gamma}_{\lambda}(b_{2}, b_{3})
-\chi^{\gamma}_{\lambda}(b_{1}\oo{\lambda}b_{2},\; b_{3})
+\chi^{\gamma}_{\lambda}(b_{1},\; b_{2}\oo{\mu}b_{3})
-\chi^{\gamma}_{\lambda}(b_{1}, b_{2})\brr{\lambda+\mu}{\gamma}b_{3}\\
=&\gamma(b_{1})\oo{\lambda}\Big(\gamma(b_{2})\oo{\mu}\gamma(b_{3})
-\gamma(b_{2}\oo{\mu} b_{3})\Big)
-\gamma(b_{1}\oo{\lambda}b_{2})\oo{\lambda+\mu}\gamma(b_{3})
+\gamma((b_{1}\oo{\lambda}b_{2})\oo{\lambda+\mu}b_{3})\\
&\ \ +\gamma(b_{1})\oo{\lambda}\gamma(b_{2}\oo{\mu}b_{3})
-\gamma(b_{1}\oo{\lambda}(b_{2}\oo{\mu}b_{3}))
-\Big(\gamma(b_{1})\oo{\lambda}\gamma(b_{2})-\gamma(b_{1}\oo{\lambda} b_{2})\Big)
\oo{\lambda+\mu}\gamma(b_{3})\\
=&0.
\end{align*}
Thus the triplet $(\blacktriangleright^{\gamma}, \blacktriangleleft^{\gamma}, \chi^{\gamma})$
is a non-abelian 2-cocycle on $B$ with values in $A$.
\end{proof}

For a split sequence in category of $\cb[\partial]$-modules, the section is not unique
in general. For two different sections, we have the following lemma.

\begin{lem} \label{lem: 2-cyc-equ}
Let $\mathcal{E}: \xymatrix@C=0.5cm{0\ar[r]& A\ar[r]^{\alpha}&E\ar[r]^{\beta}&B\ar[r]&0}$
be a non-abelian extension of associative conformal algebra $B$ by associative conformal
algebra $A$. For any two sections $\gamma$ and $\gamma'$ of $\mathcal{E}$, we have
$(\blacktriangleright^{\gamma}, \blacktriangleleft^{\gamma},
\chi^{\gamma})\approx(\blacktriangleright^{\gamma'},
\blacktriangleleft^{\gamma'}, \chi^{\gamma'})$.
\end{lem}

\begin{proof}
For the sections $\gamma$ and $\gamma'$, we define $\delta=\gamma'-\gamma:
B\rightarrow A$. Then for any $a\in A$ and $b\in B$, we have
\begin{equation*}
b\bll{\lambda}{\gamma}a=\gamma(b)\oo{\lambda}a=\gamma'(b)\oo{\lambda}a
+(\gamma-\gamma')(b)\oo{\lambda}a=b\bll{\lambda}{\gamma'}a-\delta(b)\oo{\lambda}a.
\end{equation*}
Similarly, one can check that the equations (\ref{coh7}) and (\ref{coh8}) are true.
Thus, $(\blacktriangleright^{\gamma}, \blacktriangleleft^{\gamma}, \chi^{\gamma})$ and
$(\blacktriangleright^{\gamma'}, \blacktriangleleft^{\gamma'}, \chi^{\gamma'})$ are equivalent.
\end{proof}

Therefore, we get a map from the set of all the non-abelian extensions of associative
conformal algebra $B$ by associative conformal algebra $A$ to $\HH^{2}_{nab}(B, A)$.
Next we show that this map keeps the equivalence relation.

\begin{lem} \label{lem: equ}
Let $A$ and $B$ be two associative conformal algebras. If $\mathcal{E}$ and $\mathcal{E}'$ are
two equivalent non-abelian extensions of $B$ by $A$, $\gamma$ and $\gamma'$ are the sections
of $\mathcal{E}$ and $\mathcal{E}'$ respectively. With the above notations, we get
$(\blacktriangleright^{\gamma}, \blacktriangleleft^{\gamma},
\chi^{\gamma})\approx(\blacktriangleright^{\gamma'},
\blacktriangleleft^{\gamma'}, \chi^{\gamma'})$.
\end{lem}

\begin{proof}
Since $\mathcal{E}$ and $\mathcal{E}'$ are two equivalent non-abelian extensions of
$B$ by $A$, there exists a homomorphism of associative conformal algebras
$\theta: E\rightarrow E'$ and a commutative diagram:
$$
\xymatrix{0\ar[r]& A\ar[r]^{\alpha}\ar@{=}[d]&E\ar[r]^{\beta}\ar[d]^{\theta}
&B\ar[r]\ar@{=}[d]&0\;\\ 0\ar[r]& A\ar[r]^{\alpha'}&E'\ar[r]^{\beta'}&B\ar[r]&0.}
$$
Consider the sections $\gamma: B\rightarrow E$ and $\gamma': B\rightarrow E'$, we define
$\delta=\gamma-\theta^{-1}\circ\gamma': B\rightarrow A$. Then for any $a\in A$ and
$b\in B$,
\begin{equation*}
b\bll{\lambda}{\gamma'}a=\gamma'(b)\oo{\lambda}a=\gamma(b)\oo{\lambda}a+
\theta^{-1}\circ\gamma'(b)\oo{\lambda}a-\gamma(b)\oo{\lambda}a
=b\bll{\lambda}{\gamma}a-\delta(b)\oo{\lambda}a.
\end{equation*}
Similarly, one can check the equations (\ref{coh7})-(\ref{coh8}) are true.
Thus $(\blacktriangleright^{\gamma}, \blacktriangleleft^{\gamma},
\chi^{\gamma})\approx(\blacktriangleright^{\gamma'},
\blacktriangleleft^{\gamma'}, \chi^{\gamma'})$.
\end{proof}

Let $A$ and $B$ be two associative conformal algebras. Denote by $\Ext_{nab}(B, A)$ the set of
all the equivalence classes of non-abelian extensions of $B$ by $A$, by $[\mathcal{E}]$
the equivalence class of $\mathcal{E}$ for each extension $\mathcal{E}$.
According to the above discussion, we get a map
\begin{align*}
\Phi:\quad \Ext_{nab}(B, A)\quad &\longrightarrow\quad \HH^2_{nab}(B, A)\\
[\mathcal{E}]\qquad &\;\mapsto\quad [(\blacktriangleright^{\gamma},
\blacktriangleleft^{\gamma}, \chi^{\gamma})].
\end{align*}

Conversely, for any non-abelian 2-cocycle $(\blacktriangleright, \blacktriangleleft, \chi)$
on $B$ with values in $A$, we define a conformal multiplication on $\cb[\partial]$-module
$A\oplus B$ by
$$
(a_{1}, b_{1})\oo{\lambda}(a_{2}, b_{2})=\Big(a_{1}\oo{\lambda}a_{2}+b_{1}\obl{\lambda}a_{2}
+a_{1}\obr{\lambda}b_{2}+\chi_{\lambda}(b_{1}, b_{2}),\; b_{1}\oo{\lambda}b_{2}\Big),
$$
for any $(a_{1}, b_{1})$, $(a_{2}, b_{2})\in A\oplus B$. Then we get an associative
conformal algebra, denote by $A\oplus_{(\blacktriangleright, \blacktriangleleft, \chi)}B$.
By this associative conformal algebra, we obtain a non-abelian extension of $B$ by $A$ as
follows:
$$
\mathcal{E}_{(\blacktriangleright, \blacktriangleleft, \chi)}:\qquad \xymatrix{0\ar[r]&
A\ar[r]^{\alpha\qquad}&A\oplus_{(\blacktriangleright, \blacktriangleleft, \chi)}
B\ar[r]^{\qquad\beta}&B\ar[r]&0,}
$$
where $\alpha(a)=(a, 0)$, $\beta(a, b)=b$. Moreover, we have the following lemma.

\begin{lem} \label{lem: equiv}
Let $A$ and $B$ be two associative conformal algebras. Then two non-abelian 2-cocycles
$(\blacktriangleright, \blacktriangleleft, \chi)$ and $(\bar{\blacktriangleright},
\bar{\blacktriangleleft}, \bar{\chi})$ of $B$ by $A$ are equivalent if and only if the
corresponding non-abelian extensions $\mathcal{E}_{(\blacktriangleright,
\blacktriangleleft, \chi)}$ and $\mathcal{E}_{(\bar{\blacktriangleright},
\bar{\blacktriangleleft}, \bar{\chi})}$ are equivalent.
\end{lem}

\begin{proof}
If two non-abelian 2-cocycles $(\blacktriangleright,
\blacktriangleleft, \chi)$ and $(\bar{\blacktriangleright}, \bar{\blacktriangleleft},
\bar{\chi})$ are equivalent, there is a linear map $\delta: B \rightarrow A$ such that the
equations (\ref{coh6})-(\ref{coh8}) hold. Define $\theta: A\oplus_{(\blacktriangleright,
\blacktriangleleft, \chi)}B\rightarrow A\oplus_{(\bar{\blacktriangleright},
\bar{\blacktriangleleft}, \bar{\chi})}B$ by
$$
\theta(a, b)=(a-\delta(b),\; b).
$$
It is easy to see that $\theta$ is a linear bijection with inverse $\theta^{-1}(a, b)=
(a+\delta(b),\; b)$. For any $(a_{1}, b_{1})$, $(a_{2}, b_{2})\in A\oplus_{(\blacktriangleright,
\blacktriangleleft, \chi)}B$, by equations (\ref{coh6})-(\ref{coh8}), we have
\begin{align*}
\theta((a_{1}, b_{1})\oo{\lambda}(a_{2}, b_{2}))&=\theta\Big(a_{1}\oo{\lambda}a_{2}
+b_{1}\obl{\lambda}a_{2}+a_{1}\obr{\lambda}b_{2}+\chi_{\lambda}(b_{1}, b_{2}),\;
b_{1}\oo{\lambda}b_{2}\Big)\\
&=\Big(a_{1}\oo{\lambda}a_{2}+b_{1}\obl{\lambda}a_{2}+a_{1}\obr{\lambda}b_{2}
+\chi_{\lambda}(b_{1}, b_{2})-\delta(b_{1}\oo{\lambda}b_{2}),\; b_{1}\oo{\lambda}b_{2}\Big)\\
&=\Big((a_{1}-\delta(b_{1}))\oo{\lambda}(a_{2}-\delta(b_{2}))
+b_{1}\obll{\lambda}(a_{2}-\delta(b_{2}))\\
&\qquad\qquad\qquad\qquad+(a_{1}-\delta(b_{1}))\obrr{\lambda}b_{2}
+\bar{\chi}_{\lambda}(b_{1}, b_{2}),\; b_{1}\oo{\lambda}b_{2}\Big)\\
&=\theta(a_{1}, b_{1})\oo{\lambda}\theta(a_{2}, b_{2}).
\end{align*}
Thus $\theta$ is an isomorphism of associative conformal algebras.
Moreover, one can check that the diagram:
$$
(\spadesuit)\qquad\qquad\qquad
\xymatrix{0\ar[r]& A\ar[r]^{\alpha\qquad}\ar@{=}[d]&A\oplus_{(\blacktriangleright,
\blacktriangleleft, \chi)}B\ar[r]^{\qquad\beta}\ar[d]^{\theta}
&B\ar[r]\ar@{=}[d]&0\;\\ 0\ar[r]& A\ar[r]^{\bar{\alpha}\qquad}
&A\oplus_{(\bar{\blacktriangleright},\bar{\blacktriangleleft}, \bar{\chi})}
B\ar[r]^{\qquad\bar{\beta}}&B\ar[r]&0}\qquad\qquad\quad
$$
is commutative. This means that $\mathcal{E}_{(\blacktriangleright, \blacktriangleleft, \chi)}$
and $\mathcal{E}_{(\bar{\blacktriangleright}, \bar{\blacktriangleleft}, \bar{\chi})}$ are
equivalent.

Conversely, if $\mathcal{E}_{(\blacktriangleright, \blacktriangleleft, \chi)}$
and $\mathcal{E}_{(\bar{\blacktriangleright}, \bar{\blacktriangleleft}, \bar{\chi})}$
are equivalent, by the calculation above, it is not hard to see that there is a
linear map $\theta$ such that the diagram $(\spadesuit)$ is commutative if and only if
there exists a uniquely determined linear map $\delta: B\rightarrow A$ such that
$\theta(a, b)=(a-\delta(b), b)$ for any $a\in A$ and $b\in B$.
And $\theta$ is a morphism of associative conformal algebras
if and only if $\delta$ satisfies equations (\ref{coh6})-(\ref{coh8}).
Thus $(\blacktriangleright, \blacktriangleleft, \chi)$ and $(\bar{\blacktriangleright},
\bar{\blacktriangleleft}, \bar{\chi})$ are equivalent.
\end{proof}

Thus we get a map
\begin{align*}
\Psi:\quad\HH^2_{nab}(B, A)\quad &\longrightarrow\quad\Ext_{nab}(B, A)\\
[(\blacktriangleright, \blacktriangleleft, \chi)]\quad &\;\mapsto
\qquad[\mathcal{E}_{(\blacktriangleright, \blacktriangleleft, \chi)}].
\end{align*}
And following from Lemma \ref{lem: equiv}, $\Psi$ is injective.
Next, we will show that the two mappings $\Phi$ and $\Psi$ are mutually inverse, so as to get
the main result of this section.

\begin{thm} \label{thm1}
Let $A$ and $B$ be two associative conformal algebras. There is a 1-1 correspondence between
classes of non-abelian extensions of $B$ by $A$ and elements of $\HH^{2}_{nab}(B, A)$.
In other words, $\HH^{2}_{nab}(B, A)$ classifies non-abelian extensions of $B$ by $A$.
\end{thm}

\begin{proof}
Given a non-abelian extension $\mathcal{E}: \xymatrix@C=0.5cm{0\ar[r]& A\ar[r]^{\alpha}&
E\ar[r]^{\beta}&B\ar[r]&0}$ of $B$ by $A$ and a section $\gamma$ of $\mathcal{E}$,
we get a non-abelian 2-cocycle $(\blacktriangleright^{\gamma}, \blacktriangleleft^{\gamma},
\chi^{\gamma})$. By the 2-cocycle, we obtain a non-abelian extension
$\mathcal{E}_{(\blacktriangleright^{\gamma}, \blacktriangleleft^{\gamma}, \chi^{\gamma})}$.
We define $\theta: A\oplus_{(\blacktriangleright^{\gamma},\blacktriangleleft^{\gamma},
\chi^{\gamma})}B\rightarrow E$ by
$$
\theta(a, b)=a+\gamma(b)
$$
for $a\in A$ and $b\in B$. Then $\theta$ is a bijection and for any $(a_{1}, b_{1})$,
$(a_{2}, b_{2})\in A\oplus_{(\blacktriangleright^{\gamma},\blacktriangleleft^{\gamma},
\chi^{\gamma})}B$,
\begin{align*}
\theta((a_{1}, b_{1})\oo{\lambda}(a_{2}, b_{2}))&=\theta\Big(a_{1}\oo{\lambda}a_{2}
+b_{1}\bll{\lambda}{\gamma}a_{2}+a_{1}\brr{\lambda}{\gamma}b_{2}
+\chi^{\gamma}_{\lambda}(b_{1}, b_{2}),\; b_{1}\oo{\lambda}b_{2}\Big)\\
&=a_{1}\oo{\lambda}a_{2}+\gamma(b_{1})\oo{\lambda}a_{2}+a_{1}\oo{\lambda}\gamma(b_{2})
+\gamma(b_{1})\oo{\lambda}\gamma(b_{2})\\
&=\theta(a_{1}, b_{1})\oo{\lambda}\theta(a_{2}, b_{2}).
\end{align*}
Moreover, one can check that the diagram:
$$
\xymatrix{0\ar[r]& A\ar[r]^{\alpha'\qquad}\ar@{=}[d]&
A\oplus_{(\blacktriangleright^{\gamma},\blacktriangleleft^{\gamma}, \chi^{\gamma})}B
\ar[r]^{\qquad\beta'}\ar[d]^{\theta}&B\ar[r]\ar@{=}[d]&0\\
0\ar[r]& A\ar[r]^{\alpha}&E\ar[r]^{\beta} &B\ar[r]&0 }
$$
is commutative. This means that $\Psi\circ\Phi(\mathcal{E})$ and $\mathcal{E}$ are equivalent.
Thus $\Psi$ is surjective, and so that it is a bijection.
\end{proof}

At the end of this section, we apply the conclusion of Theorem \ref{thm1} to
the $\cb[\partial]$-split abelian extension of associative conformal algebras.
If $B$ as $\cb[\partial]$-module is projective and the multiplication in $A$
is trivial, the Definition \ref{def:coho} is just the definition of $\cb[\partial]$-split
abelian extension of $B$ by $A$ in \cite{Do}. In this case, for an extension
$\mathcal{E}: \xymatrix@C=0.5cm{0\ar[r]& A\ar[r]&E\ar[r]&B\ar[r]&0}$,
the corresponding non-abelian 2-cocycle $(\blacktriangleright, \blacktriangleleft, \chi)$
induces a bimodule structure on $A$ by equations (\ref{coh1})-(\ref{coh3}), and gives
an element $\chi\in\mathcal{Z}^{2}(B, A)$ by equation (\ref{coh5}).
The equation (\ref{coh8}) means that two extensions $\mathcal{E}$ and $\bar{\mathcal{E}}$
are equivalent if and only if $\chi-\bar{\chi}\in\mathcal{B}^{2}(B, A)$. Thus we obtain
the following corollary.

\begin{cor}[\cite{Do}] \label{cor1}
Let $B$ be an associative conformal algebra, which is projective as $\cb[\partial]$-module.
For any bimodule $A$ over $B$, the Hochschild cohomology $\HH^{2}(B, A)$ classifies
abelian extensions of $B$ by $A$.
\end{cor}

%%%%%%%%%%%%%%%%%%%%%%%%%%%%%%%%%%%%%%%%%%%%%%%%%%%%%%%%%%%%%%%%%%%%%%%%%%%%%%%%%
%    section  5   Non-abelian extensions in terms of Deligne groupoid
%%%%%%%%%%%%%%%%%%%%%%%%%%%%%%%%%%%%%%%%%%%%%%%%%%%%%%%%%%%%%%%%%%%%%%%%%%%%%%%%%%%%%%

\section{Non-abelian extensions in terms of Deligne groupoid}\label{sec:Deli}

In this section, we construct a subalgebra $\mathcal{L}$ of differential graded algebra
$(\mathcal{C}^{\bullet+1}(A\oplus B, A\oplus B),\; [-,-],\; d_{\ast})$ for two
associative conformal algebras $A$ and $B$. We show that the non-abelian cocycles are
in bijection with the Maurer-Cartan elements of $\mathcal{L}$, and get that the equivalence
relation on $\mathcal{Z}^{2}_{nab}(B, A)$ can be interpreted as gauge equivalence relation on
${\rm MC}(\mathcal{L})$.

Let $\mathcal{E}: \xymatrix@C=0.5cm{0\ar[r]& A\ar[r]^{\alpha}&E\ar[r]^{\beta}&B\ar[r]&0}$
be a non-abelian extension of associative conformal algebra $B$ by associative conformal
algebra $A$. Suppose $\gamma$ is a section of $\mathcal{E}$.
Then in the category of $\cb[\partial]$-modules, we have a commutative diagram:
$$
\xymatrix{0\ar[r]& A\ar[r]^{\alpha}\ar@{=}[d]&E\ar[r]^{\beta}\ar@{=}[d]
&B\ar[r]\ar[d]^{\cong}&0\;\\
0\ar[r]& \bar{A}\ar[r]^{\iota_{\bar{A}}\quad}&\bar{A}\oplus\bar{B}
\ar[r]^{\quad\rho_{\bar{B}}}&\bar{B}\ar[r]&0,}
$$
where $\iota_{\bar{A}}(a)=(a, 0)$, $\rho_{\bar{B}}(a, b)=b$,
$\bar{A}$ is the $\cb[\partial]$-module of image of $A$ under $\alpha$,
$\bar{B}$ is the $\cb[\partial]$-module of image of $B$ under $\gamma$.
Using these isomorphisms, we can transfer the associative conformal algebra structure
of $A$, $B$ and $E$ to $\bar{A}$, $\bar{B}$ and $\bar{A}\oplus\bar{B}$, such that this
commutative graph holds in category of associative conformal algebras.
Denote the multiplication in $\bar{A}\oplus\bar{B}$ (or in $E$) by $\mathfrak{M}$, i.e.,
$\mathfrak{M}_{\lambda}((a_{1}, b_{1}), (a_{2}, b_{2}))=(a_{1}, b_{1})\oo{\lambda}
(a_{2}, b_{2})$. Then we have $\mathfrak{M}=\mathfrak{M}_{\bar{A}\bar{A}}^{\bar{A}}
+\mathfrak{M}_{\bar{A}\bar{B}}^{\bar{A}}+\mathfrak{M}_{\bar{B}\bar{A}}^{\bar{A}}
+\mathfrak{M}_{\bar{B}\bar{B}}^{\bar{A}}+\mathfrak{M}_{\bar{A}\bar{A}}^{\bar{B}}
+\mathfrak{M}_{\bar{A}\bar{B}}^{\bar{B}}+\mathfrak{M}_{\bar{B}\bar{A}}^{\bar{B}}
+\mathfrak{M}_{\bar{B}\bar{B}}^{\bar{B}}$, where
$(\mathfrak{M}_{\bar{A}\bar{B}}^{\bar{A}})_{\lambda}
((a_{1}, b_{1}), (a_{2}, b_{2}))=\rho_{\bar{A}}((a_{1}, 0)\oo{\lambda}(0, b_{2}))$,
and others can be defined similarly. For these components, we have the following lemma.

\begin{lem}\label{lem: compo}
With the above notations, we have the following conclusions.

$(i)$ One can identifies $\mathfrak{M}_{\bar{B}\bar{B}}^{\bar{B}}$
with multiplication in $B$, identifies $\mathfrak{M}_{\bar{A}\bar{A}}^{\bar{A}}$ with
multiplication in $A$, and $\mathfrak{M}_{\bar{A}\bar{B}}^{\bar{B}}
=\mathfrak{M}_{\bar{B}\bar{A}}^{\bar{B}}=\mathfrak{M}_{\bar{A}\bar{A}}^{\bar{B}}=0$.

$(ii)$ Denote by $(\blacktriangleright, \blacktriangleleft, \chi)$ the non-abelian 2-cocycle
on $B$ with values in $A$ corresponding to $\mathcal{E}$. Then one can identify
$\mathfrak{M}_{\bar{B}\bar{A}}^{\bar{A}}$, $\mathfrak{M}_{\bar{A}\bar{B}}^{\bar{A}}$ and
$\mathfrak{M}_{\bar{B}\bar{B}}^{\bar{A}}$ with $\blacktriangleright$, $\blacktriangleleft$
and $\chi$ respectively.
\end{lem}

\begin{proof}
$(i)$ Since $\beta\circ\gamma=\id_{B}$, we can identify
$\mathfrak{M}_{\bar{B}\bar{B}}^{\bar{B}}$ with multiplication in $B$.
And similarly for $\mathfrak{M}_{\bar{A}\bar{A}}^{\bar{A}}$ and multiplication in $A$.
Moreover, since $\rho_{\bar{B}}$ is a morphism of associative conformal algebras and
$\rho_{\bar{B}}(a, 0)=0$, for any $(a_{1}, b_{1})$, $(a_{2}, b_{2})\in \bar{A}\oplus \bar{B}$,
$$
(\mathfrak{M}_{\bar{B}\bar{A}}^{\bar{B}})_{\lambda}((a_{1}, b_{1}),\; (a_{2}, b_{2}))
=\rho_{\bar{B}}\circ\mathfrak{M}_{\lambda}((0, b_{1}),\; (a_{2}, 0))
=\rho_{\bar{B}}(0, b_{1})\oo{\lambda}\rho_{\bar{B}}(a_{2}, 0)=0.
$$
Similarly, $\mathfrak{M}_{\bar{A}\bar{B}}^{\bar{B}}=\mathfrak{M}_{\bar{A}\bar{A}}^{\bar{B}}=0$.

$(ii)$ If $(\blacktriangleright, \blacktriangleleft, \chi)$ is the non-abelian 2-cocycle
on $B$ with values in $A$ corresponding to $\mathcal{E}$. Then the multiplication in
$\bar{A}\oplus\bar{B}$ (or in $E$) is given by
$$
\mathfrak{M}_{\lambda}((a_{1}, b_{1}),\; (a_{2}, b_{2}))=\Big(a_{1}\oo{\lambda}a_{2}
+b_{1}\obl{\lambda}a_{2}+a_{1}\obr{\lambda}b_{2}+\chi_{\lambda}(b_{1}, b_{2}),\; b_{1}
\oo{\lambda}b_{2}\Big),
$$
for any $a_{1}, a_{2}\in A$, $b_{1}, b_{2}\in B$. Thus,
$$
(\mathfrak{M}_{\bar{B}\bar{A}}^{\bar{A}})_{\lambda}((0, b_{1}),\; (a_{2}, 0))
=\rho_{\bar{A}}\circ\mathfrak{M}_{\lambda}((0, b_{1}),\; (a_{2}, 0))
=b_{1}\obl{\lambda}a_{2}.
$$
Hence we can identify $\mathfrak{M}_{\bar{B}\bar{A}}^{\bar{A}}$ with $\blacktriangleright$.
Similarly, we identify $\mathfrak{M}_{\bar{A}\bar{B}}^{\bar{A}}$ with
$\blacktriangleleft$, $\mathfrak{M}_{\bar{B}\bar{B}}^{\bar{A}}$ with $\chi$.
\end{proof}

Next, we consider the associator $\mathcal{A}ss$ of the multiplication $\mathfrak{M}$ in
$\bar{A}\oplus\bar{B}$, i.e., $\mathcal{A}ss=\mathfrak{M}\otimes\id-\id\otimes\mathfrak{M}$.
Consider the components of $\mathcal{A}ss$, we have the following lemma.

\begin{pro}\label{pro: ass}
With the above notations, we have
\begin{itemize}
\item[(i)] $\mathcal{A}ss^{\bar{B}}_{\bar{B}\bar{B}\bar{B}}=0$ is equivalent to the
           multiplication in $B$ is associative;
\item[(ii)] $\mathcal{A}ss^{\bar{A}}_{\bar{A}\bar{A}\bar{A}}=0$ is equivalent to the
           multiplication in $A$ is associative;
\item[(iii)] $\mathcal{A}ss^{\bar{A}}_{\bar{B}\bar{B}\bar{A}}=0$ is equivalent to the
           equation (\ref{coh1});
\item[(iv)] $\mathcal{A}ss^{\bar{A}}_{\bar{A}\bar{B}\bar{B}}=0$ is equivalent to the
           equation (\ref{coh2});
\item[(v)] $\mathcal{A}ss^{\bar{A}}_{\bar{B}\bar{A}\bar{B}}=0$ is equivalent to the
           equation (\ref{coh3});
\item[(vi)] $\mathcal{A}ss^{\bar{A}}_{\bar{A}\bar{A}\bar{B}}=0$ is equivalent to the
           equation (\ref{coh4});
\item[(vii)] $\mathcal{A}ss^{\bar{A}}_{\bar{B}\bar{A}\bar{A}}=0$ is equivalent to the
           equation (\ref{coh4'});
\item[(viii)] $\mathcal{A}ss^{\bar{A}}_{\bar{A}\bar{B}\bar{A}}=0$ is equivalent to the
           equation (\ref{coh4''});
\item[(ix)] $\mathcal{A}ss^{\bar{A}}_{\bar{B}\bar{B}\bar{B}}=0$ is equivalent to the
           equation (\ref{coh5}).
\end{itemize}
\end{pro}

\begin{proof}
Using the conclusion of Lemma \ref{lem: compo}, this proposition can be obtained by direct
calculation. The details are as follows. Let $a_{1}, a_{2}, a_{3}\in A$ and
$b_{1}, b_{2}, b_{3}\in B$.

$(i)$ Note that
\begin{align*}
&(\mathcal{A}ss^{\bar{B}}_{\bar{B}\bar{B}\bar{B}})_{\lambda, \mu}\Big((a_{1}, b_{1}),\;
(a_{2}, b_{2}),\; (a_{3}, b_{3})\Big)\\
=&(\mathfrak{M}_{\bar{B}\bar{B}}^{\bar{B}})_{\lambda+\mu}
\Big((\mathfrak{M}_{\bar{B}\bar{B}}^{\bar{B}})_{\lambda}
((a_{1}, b_{1}),\; (a_{2}, b_{2})),\; (a_{3}, b_{3})\Big)
+(\mathfrak{M}_{\bar{A}\bar{B}}^{\bar{B}})_{\lambda+\mu}
(\Big(\mathfrak{M}_{\bar{B}\bar{B}}^{\bar{A}})_{\lambda}
((a_{1}, b_{1}),\; (a_{2}, b_{2})),\; (a_{3}, b_{3})\Big)\\
&-(\mathfrak{M}_{\bar{B}\bar{B}}^{\bar{B}})_{\lambda}\Big((a_{1}, b_{1}),\;
(\mathfrak{M}_{\bar{B}\bar{B}}^{\bar{B}})_{\mu}((a_{2}, b_{2}),\; (a_{3}, b_{3}))\Big)
-(\mathfrak{M}_{\bar{B}\bar{A}}^{\bar{B}})_{\lambda}\Big((a_{1}, b_{1}),\;
(\mathfrak{M}_{\bar{B}\bar{B}}^{\bar{A}})_{\mu}((a_{2}, b_{2}),\; (a_{3}, b_{3}))\Big)\\
=&(b_{1}\oo{\lambda}b_{2})\oo{\lambda+\mu}b_{3}-b_{1}\oo{\lambda}(b_{2}\oo{\mu}b_{3}).
\end{align*}
Thus $\mathcal{A}ss^{\bar{B}}_{\bar{B}\bar{B}\bar{B}}=0$ is equivalent to the multiplication
in $B$ is associative. Similarly, we get $(ii)$.

$(iii)$ Since
\begin{align*}
&(\mathcal{A}ss^{\bar{A}}_{\bar{B}\bar{B}\bar{A}})_{\lambda, \mu}
\Big((a_{1}, b_{1}),\; (a_{2}, b_{2}),\; (a_{3}, b_{3})\Big)\\
=&(\mathfrak{M}_{\bar{B}\bar{A}}^{\bar{A}})_{\lambda+\mu}
\Big((\mathfrak{M}_{\bar{B}\bar{B}}^{\bar{B}})_{\lambda}
((a_{1}, b_{1}),\; (a_{2}, b_{2})),\; (a_{3}, b_{3})\Big)
+(\mathfrak{M}_{\bar{A}\bar{A}}^{\bar{A}})_{\lambda+\mu}
\Big((\mathfrak{M}_{\bar{B}\bar{B}}^{\bar{A}})_{\lambda}
((a_{1}, b_{1}),\; (a_{2}, b_{2})),\; (a_{3}, b_{3})\Big)\\
&-(\mathfrak{M}_{\bar{B}\bar{B}}^{\bar{A}})_{\lambda}\Big((a_{1}, b_{1}),\;
(\mathfrak{M}_{\bar{B}\bar{A}}^{\bar{B}})_{\mu}((a_{2}, b_{2}),\; (a_{3}, b_{3}))\Big)
-(\mathfrak{M}_{\bar{B}\bar{A}}^{\bar{A}})_{\lambda}\Big((a_{1}, b_{1}),\;
(\mathfrak{M}_{\bar{B}\bar{A}}^{\bar{A}})_{\mu}((a_{2}, b_{2}),\; (a_{3}, b_{3}))\Big)\\
=&(b_{1}\oo{\lambda}b_{2})\obl{\lambda+\mu}a_{3}
+\chi_{\lambda}(b_{1}, b_{2})\oo{\lambda+\mu}a_{3}
-b_{1}\obl{\lambda}(b_{2}\obl{\mu}a_{3}).
\end{align*}
That is to say, $\mathcal{A}ss^{\bar{B}}_{\bar{B}\bar{B}\bar{B}}=0$ is equivalent to
equation (\ref{coh1}). Similarly, we can obtain $(iv)$-$(ix)$.
\end{proof}

\begin{rmk}\label{rmk2}
By direct calculation, we have $\mathcal{A}ss^{\bar{B}}_{\bar{A}\bar{B}\bar{B}}
=\mathcal{A}ss^{\bar{B}}_{\bar{A}\bar{A}\bar{B}}
=\mathcal{A}ss^{\bar{B}}_{\bar{A}\bar{B}\bar{A}}
=\mathcal{A}ss^{\bar{B}}_{\bar{B}\bar{A}\bar{A}}
=\mathcal{A}ss^{\bar{B}}_{\bar{A}\bar{A}\bar{A}}=0$.
\end{rmk}

Given two associative conformal algebras $A$ and $B$, we get an associative conformal
algebra structure on the direct sum of $\cb[\partial]$-module $A\oplus B$ by the
conformal multiplication $\mathfrak{m}_{A\oplus B}$:
\begin{align*}
(\mathfrak{m}_{A\oplus B})_{\lambda}\Big((a_{1}, b_{1}),\; (a_{2}, b_{2})\Big)
&=(a_{1}, b_{1})\oo{\lambda}(a_{2}, b_{2})\\
&=\Big((\mathfrak{m}_{A})_{\lambda}(a_{1}, a_{2}),\;
(\mathfrak{m}_{B})_{\lambda}(b_{1}, b_{1})\Big)
=\Big(a_{1}\oo{\lambda}a_{2},\; b_{1}\oo{\lambda}b_{2}\Big),
\end{align*}
for $a_{1}, a_{2}\in A$ and $b_{1}, b_{2}\in B$, where $\mathfrak{m}_{A}$ and
$\mathfrak{m}_{B}$ are the conformal multiplications of $A$ and $B$ respectively.
Considering the Hochschild complex of $A\oplus B$,
we get a differential graded Lie algebra $(\mathcal{C}^{\bullet+1}(A\oplus B, A\oplus B),\;
[-,-],\; \bar{d}_{\ast})$, where $[-,-]$ is the Gerstenhaber bracket, $\bar{d}_{i}=
[\mathfrak{m}_{A\oplus B}, -]=[\mathfrak{m}_{A}+\mathfrak{m}_{B}, -]$.
View $A$ as an $A\oplus B$-module via the action of $A$ on itself, and consider the Hochschild
complex of $A\oplus B$ with coefficients in $A$, we define
$$
\mathcal{L}=\bigoplus_{i\geq0}\mathcal{L}^{i},\qquad\qquad
\mathcal{L}^{i}=\bigoplus_{m+n=i+1, m\geq0, n>0}\mathcal{L}^{m,n},
$$
where $\mathcal{L}^{m,n}$ is the set all conformal sesquilinear maps from
$(A^{\otimes m}\otimes B^{\otimes n})\oplus(A^{\otimes(m-1)}\otimes B\otimes A\oplus
B^{\otimes(n-1)})\oplus\cdots\oplus (B^{\otimes n}\oplus A^{\otimes m})$ to $A$.
Then clearly, $\mathcal{L}^{i}\subseteq\mathcal{C}^{i+1}(A\oplus B, A)
\subseteq\mathcal{C}^{i+1}(A\oplus B, A\oplus B)$ and $\mathcal{C}^{i+1}(A\oplus B, A)=
\mathcal{L}^{i}\oplus\mathcal{C}^{i+1}(A, A)$. We also denote the restriction of
Gerstenhaber bracket and differential on $\mathcal{L}$ by $[-, -]$ and $\bar{d}_{\ast}$.

\begin{pro}\label{pro: DGlie}
With the above notations, the triples $(\mathcal{L},\; [-,-],\; \bar{d}_{\ast})$ is a
sub-differential graded Lie algebra of $(\mathcal{C}^{\bullet+1}(A\oplus B, A\oplus B),\;
[-,-],\; \bar{d}_{\ast})$. Its degree 0 part is abelian.
\end{pro}

\begin{proof}
$(i)$ $\mathcal{L}$ is closed under the Gerstenhaber bracket.
Denote by $L_{l,k}$ the set all conformal sesquilinear maps from
$X_{1}\otimes X_{2}\otimes\cdots\otimes X_{n}$ to $A$, where $A$ appears $l$-times and
$B$ appears $k$-times in $X_{1}, X_{2},\cdots X_{n}$ and $l+k=n$.
Then any $f\in\mathcal{L}^{m,n}$ can be decomposed as the sum of $f_{i,j}$,
where $f_{i,j}\in L_{ij}$, $i+j=m+n$, $j>0$.
For any $f\in\mathcal{L}^{m,n}$ and $g\in\mathcal{L}^{m',n'}$, we get
$$
[f, g]\in\mathcal{L}^{m+m'-1,n+n'},
$$
since $[f_{i,j}, g_{i',j'}]\in L_{i+i'-1, j+j'}$ by the definition of the Gerstenhaber bracket.
Thus, $\mathcal{L}$ is closed under the Gerstenhaber bracket.

$(ii)$ $\mathcal{L}$ is closed under the differential. For any $f\in\mathcal{L}^{m,n}\subset
\mathcal{L}^{i}$, $\bar{d}_{i+1}(f)=[\mathfrak{m}_{A}+\mathfrak{m}_{B}, f]$. Note that
$f$ can be decomposed as the sum of $f_{i,j}$, $i+j=m+n$, $j>0$, and
$$
[\mathfrak{m}_{A}, f_{ij}]\in L_{i+1, j},\qquad\qquad
[\mathfrak{m}_{B}, f_{ij}]\in L_{i, j+1},
$$
$L_{i+1, j}\subseteq\mathcal{L}^{i+1}$ and $L_{i, j+1}\subseteq\mathcal{L}^{i+1}$, we get
$\bar{d}_{i+1}(f_{ij})\in\mathcal{L}^{i+1}$, and so that $\bar{d}_{i+1}(f)\in\mathcal{L}^{i+1}$.
That is to say, $\mathcal{L}$ is closed under the differential.

Finally, note that $\mathcal{L}^{0}$ is the set of all the $\cb[\partial]$-module homomorphisms
from $B$ to $A$, we get $\mathcal{L}^{0}$ is abelian under the Gerstenhaber bracket.
\end{proof}

\begin{lem}\label{lem: MC}
Let $A$ and $B$ be two associative conformal algebras on $\cb[\partial]$-modules $\bar{A}$
and $\bar{B}$ respectively. Then we have
$$
\mathfrak{e}:=\mathcal{A}ss^{\bar{A}}_{\bar{B}\bar{B}\bar{B}}
+\mathcal{A}ss^{\bar{A}}_{\bar{B}\bar{B}\bar{A}}
+\mathcal{A}ss^{\bar{A}}_{\bar{B}\bar{A}\bar{B}}+\mathcal{A}ss^{\bar{A}}_{\bar{A}\bar{B}\bar{B}}
+\mathcal{A}ss^{\bar{A}}_{\bar{A}\bar{A}\bar{B}}+\mathcal{A}ss^{\bar{A}}_{\bar{A}\bar{B}\bar{A}}
+\mathcal{A}ss^{\bar{A}}_{\bar{B}\bar{A}\bar{A}}=0
$$
if and only if $\mathfrak{c}:=\mathfrak{M}_{\bar{B}\bar{B}}^{\bar{A}}
+\mathfrak{M}_{\bar{B}\bar{A}}^{\bar{A}}+\mathfrak{M}_{\bar{A}\bar{B}}^{\bar{A}}\in
{\rm MC}(\mathcal{L})$.
\end{lem}

\begin{proof}
For any $e_{1}=(a_{1}, b_{1})$, $e_{2}=(a_{2}, b_{2})$, $e_{3}=(a_{3}, b_{3})\in A\oplus B$,
we need compute $(\bar{d}_{2}(\mathfrak{c})+\frac{1}{2}[\mathfrak{c}, \mathfrak{c}])_{\lambda,
\mu}(e_{1}, e_{2}, e_{3})$. Note that $\bar{d}_{2}(\mathfrak{c})=[\mathfrak{m}_{A},
\mathfrak{c}]+[\mathfrak{m}_{B}, \mathfrak{c}]$,
and we can identify $\mathfrak{M}_{\bar{B}\bar{B}}^{\bar{B}}$ and
$\mathfrak{M}_{\bar{A}\bar{A}}^{\bar{A}}$ with $\mathfrak{m}_{B}$ and $\mathfrak{m}_{A}$
respectively, we have
\begin{align*}
([\mathfrak{m}_{A},\; \mathfrak{c}])_{\lambda, \mu}(e_{1}, e_{2}, e_{3})
=&\; (\mathfrak{M}_{\bar{A}\bar{A}}^{\bar{A}})_{\lambda+\mu}
\Big((\mathfrak{M}_{\bar{B}\bar{B}}^{\bar{A}})_{\lambda}(e_{1}, e_{2}),\; e_{3}\Big)
-(\mathfrak{M}_{\bar{A}\bar{A}}^{\bar{A}})_{\lambda}\Big(e_{1},\;
(\mathfrak{M}_{\bar{B}\bar{B}}^{\bar{A}})_{\mu}(e_{2}, e_{3})\Big)\\
&+(\mathfrak{M}_{\bar{A}\bar{A}}^{\bar{A}})_{\lambda+\mu}
\Big((\mathfrak{M}_{\bar{B}\bar{A}}^{\bar{A}})_{\lambda}(e_{1}, e_{2}),\; e_{3}\Big)
-(\mathfrak{M}_{\bar{A}\bar{A}}^{\bar{A}})_{\lambda}\Big(e_{1},\;
(\mathfrak{M}_{\bar{B}\bar{A}}^{\bar{A}})_{\mu}(e_{2}, e_{3})\Big)\\
&+(\mathfrak{M}_{\bar{A}\bar{A}}^{\bar{A}})_{\lambda+\mu}
\Big((\mathfrak{M}_{\bar{A}\bar{B}}^{\bar{A}})_{\lambda}(e_{1}, e_{2}),\; e_{3}\Big)
-(\mathfrak{M}_{\bar{A}\bar{A}}^{\bar{A}})_{\lambda}\Big(e_{1},\;
(\mathfrak{M}_{\bar{A}\bar{B}}^{\bar{A}})_{\mu}(e_{2}, e_{3})\Big)\\
&+(\mathfrak{M}_{\bar{A}\bar{B}}^{\bar{A}})_{\lambda+\mu}
\Big((\mathfrak{M}_{\bar{A}\bar{A}}^{\bar{A}})_{\lambda}(e_{1}, e_{2}),\; e_{3}\Big)
-(\mathfrak{M}_{\bar{B}\bar{A}}^{\bar{A}})_{\lambda}\Big(e_{1},\;
(\mathfrak{M}_{\bar{A}\bar{A}}^{\bar{A}})_{\mu}(e_{2}, e_{3})\Big),
\end{align*}
\begin{align*}
([\mathfrak{m}_{B},\; \mathfrak{c}])_{\lambda, \mu}(e_{1}, e_{2}, e_{3})
=&\; (\mathfrak{M}_{\bar{B}\bar{B}}^{\bar{A}})_{\lambda+\mu}
\Big((\mathfrak{M}_{\bar{B}\bar{B}}^{\bar{B}})_{\lambda}(e_{1}, e_{2}),\; e_{3}\Big)
-(\mathfrak{M}_{\bar{B}\bar{B}}^{\bar{A}})_{\lambda}\Big(e_{1},\;
(\mathfrak{M}_{\bar{B}\bar{B}}^{\bar{B}})_{\mu}(e_{2}, e_{3})\Big)\\
& +(\mathfrak{M}_{\bar{B}\bar{A}}^{\bar{A}})_{\lambda+\mu}
\Big((\mathfrak{M}_{\bar{B}\bar{B}}^{\bar{B}})_{\lambda}(e_{1}, e_{2}),\; e_{3}\Big)
-(\mathfrak{M}_{\bar{A}\bar{B}}^{\bar{A}})_{\lambda}\Big(e_{1},\;
(\mathfrak{M}_{\bar{B}\bar{B}}^{\bar{B}})_{\mu}(e_{2}, e_{3})\Big),
\end{align*}
and
\begin{align*}
(\mbox{$\frac{1}{2}$}[\mathfrak{c},\; \mathfrak{c}])_{\lambda, \mu}(e_{1}, e_{2}, e_{3})
=&\; (\mathfrak{M}_{\bar{A}\bar{B}}^{\bar{A}})_{\lambda+\mu}
\Big((\mathfrak{M}_{\bar{B}\bar{B}}^{\bar{A}})_{\lambda}(e_{1}, e_{2}),\; e_{3}\Big)
-(\mathfrak{M}_{\bar{B}\bar{A}}^{\bar{A}})_{\lambda}\Big(e_{1},\;
(\mathfrak{M}_{\bar{B}\bar{B}}^{\bar{A}})_{\mu}(e_{2}, e_{3})\Big)\\
& +(\mathfrak{M}_{\bar{A}\bar{B}}^{\bar{A}})_{\lambda+\mu}
\Big((\mathfrak{M}_{\bar{B}\bar{A}}^{\bar{A}})_{\lambda}(e_{1}, e_{2}),\; e_{3}\Big)
-(\mathfrak{M}_{\bar{B}\bar{A}}^{\bar{A}})_{\lambda}\Big(e_{1},\;
(\mathfrak{M}_{\bar{B}\bar{A}}^{\bar{A}})_{\mu}(e_{2}, e_{3})\Big)\\
& +(\mathfrak{M}_{\bar{A}\bar{B}}^{\bar{A}})_{\lambda+\mu}
\Big((\mathfrak{M}_{\bar{A}\bar{B}}^{\bar{A}})_{\lambda}(e_{1}, e_{2}),\; e_{3}\Big)
-(\mathfrak{M}_{\bar{B}\bar{A}}^{\bar{A}})_{\lambda}\Big(e_{1},\;
(\mathfrak{M}_{\bar{A}\bar{B}}^{\bar{A}})_{\mu}(e_{2}, e_{3})\Big).
\end{align*}
Compared with the calculation of Proposition \ref{pro: ass}, we get
$$
\Big([\mathfrak{m}_{A},\; \mathfrak{c}]+[\mathfrak{m}_{B},\; \mathfrak{c}]
+\mbox{$\frac{1}{2}$}[\mathfrak{c},\; \mathfrak{c}]\Big)_{\lambda, \mu}(e_{1}, e_{2}, e_{3})
=\mathfrak{e}_{\lambda, \mu}(e_{1}, e_{2}, e_{3}).
$$
Thus we obtain the lemma.
\end{proof}

Now we can give a 1-1 correspondence between the non-abelian 2-cocycles on $B$
with values in $A$ and Maurer-Cartan elements of $\mathcal{L}$.

\begin{pro}\label{pro: corres}
Let $A$ and $B$ be two associative conformal algebras, and $\mathcal{L}$ be the differential
graded Lie algebra defined as above. Then we have a 1-1 correspondence
\begin{align*}
\mathcal{Z}^{2}_{nab}(B, A)\quad&\longleftrightarrow\quad {\rm MC}(\mathcal{L})\\
(\blacktriangleright, \blacktriangleleft, \chi)\quad &\;\;\leftrightarrow \;\;\quad
\chi+\blacktriangleright+\blacktriangleleft.
\end{align*}
\end{pro}

\begin{proof}
By Theorem \ref{thm1}, we know that an element in $\mathcal{Z}^{2}_{nab}(B, A)$ can be view as
a non-abelian extension of $B$ by $A$. Given a non-abelian extension of $B$ by $A$,
since the extension algebra is associative, we get the
associator $\mathcal{A}ss=\mathcal{A}ss^{\bar{A}}_{\bar{A}\bar{A}\bar{A}}
+\mathfrak{e}+\mathcal{A}ss^{\bar{B}}_{\bar{B}\bar{B}\bar{B}}=0$.
Since $A$ and $B$ are associative, i.e.,
$\mathcal{A}ss^{\bar{A}}_{\bar{A}\bar{A}\bar{A}}=0$ and
$\mathcal{A}ss^{\bar{B}}_{\bar{B}\bar{B}\bar{B}}=0$, we have $\mathfrak{e}=0$.
By Lemma \ref{lem: MC}, there exists a unique Maurer-Cartan element $\mathfrak{c}=
\mathfrak{M}_{\bar{B}\bar{B}}^{\bar{A}}+\mathfrak{M}_{\bar{B}\bar{A}}^{\bar{A}}
+\mathfrak{M}_{\bar{A}\bar{B}}^{\bar{A}}=\chi+\blacktriangleright+\blacktriangleleft$ of
$\mathcal{L}$ such that $\mathfrak{e}=0$. Thus we obtain the 1-1 correspondence.
\end{proof}

Next we show that the 1-1 correspondence in the proposition above also
keeps the equivalence relations on $\mathcal{Z}^{2}_{nab}(B, A)$ and ${\rm MC}(\mathcal{L})$.
We then give the main theorem of this section.

\begin{thm}\label{thm2}
Let $A$ and $B$ be two associative conformal algebras, and $\mathcal{L}$ be the differential
graded Lie algebra defined as above. Then we have a 1-1 correspondence
$$
\HH^{2}_{nab}(B, A)\cong \mathcal{MC}(\mathcal{L}).
$$
\end{thm}

\begin{proof}
We need show that equivalence relation on non-abelian 2-cocycles coincides with
gauge relation on ${\rm MC}(\mathcal{L})$. Recall that two elements $\mathfrak{c}$ and
$\tilde{\mathfrak{c}}$ in ${\rm MC}(L)$ are equivalent if there exists $\xi\in\mathcal{L}^{0}
\subseteq\Hom_{\cb[\partial]}(A\oplus B, A)$ such that
$$
\tilde{\mathfrak{c}}=e^{\ad_{\xi}}\mathfrak{c}-\frac{e^{\ad_{\xi}}-1}{\ad_{\xi}}\bar{d}_{1}(\xi).
$$
We consider $\mathfrak{c}:=\chi+\blacktriangleright+\blacktriangleleft
=\mathfrak{M}_{\bar{B}\bar{B}}^{\bar{A}}+\mathfrak{M}_{\bar{B}\bar{A}}^{\bar{A}}
+\mathfrak{M}_{\bar{A}\bar{B}}^{\bar{A}}$. Since $\mathcal{L}^{0}$
is the set of all $\cb[\partial]$-module homomorphisms from $B$ to $A$, $\xi$ is ad-nilpotent.
Thus for any $e_{1}=(a_{1}, b_{1})$, $e_{2}=(a_{2}, b_{2})\in A\oplus B$, we have
\begin{align*}
e^{\ad_{\xi}}(\mathfrak{c})_{\lambda}(e_{1}, e_{2})&=\Big(\mathfrak{c}+[\xi, \mathfrak{c}]
+\mbox{$\frac{1}{2}$}[\xi, [\xi, \mathfrak{c}]]+\cdots\Big)_{\lambda}(e_{1}, e_{2})\\
&=b_{1}\obl{\lambda}a_{2}+a_{1}\obr{\lambda}b_{2}+\chi_{\lambda}(b_{1}, b_{2})
+([\xi, \mathfrak{c}])_{\lambda}(e_{1}, e_{2}).
\end{align*}
Note that $[\xi,\chi]=0$, since $\xi$ and $\chi$ take values in $A$ and $A[\lambda]$, and
\begin{align*}
([\xi, \blacktriangleright])_{\lambda}(e_{1}, e_{2})&=\xi(\obl{\lambda}(e_{1}, e_{2}))
-\obl{\lambda}(\xi(e_{1}), e_{2})-\obl{\lambda}(e_{1}, \xi(e_{2}))
=-b_{1}\obl{\lambda}\xi(b_{2}),\\
([\xi, \blacktriangleleft])_{\lambda}(e_{1}, e_{2})&=\xi(\obr{\lambda}(e_{1}, e_{2}))
-\obr{\lambda}(\xi(e_{1}), e_{2})-\obr{\lambda}(e_{1}, \xi(e_{2}))
=-\xi(b_{1})\obr{\lambda}b_{2},
\end{align*}
we get
$$
e^{\ad_{\xi}}(\mathfrak{c})_{\lambda}(e_{1}, e_{2})=b_{1}\obl{\lambda}a_{2}+a_{1}
\obr{\lambda}b_{2}+\chi_{\lambda}(b_{1}, b_{2})-b_{1}\obl{\lambda}\xi(b_{2})
-\xi(b_{1})\obr{\lambda}b_{2}.
$$

Next, since $\bar{d}_{1}(\xi)_{\lambda}(e_{1}, e_{2})=([\mathfrak{m}_{A}
+\mathfrak{m}_{B},\; \xi])_{\lambda}(e_{1}, e_{2})$ and
\begin{align*}
([\mathfrak{m}_{A}, \xi])_{\lambda}(e_{1}, e_{2})&=(\mathfrak{m}_{A})_{\lambda}
(\xi(e_{1}), e_{2}))+(\mathfrak{m}_{A})_{\lambda}(e_{1}, \xi(e_{2}))
-\xi((\mathfrak{m}_{A})_{\lambda}(e_{1}, e_{2}))\\
&=\xi(b_{1})\oo{\lambda}a_{2}+a_{1}\oo{\lambda}\xi(b_{2}),\\
([\mathfrak{m}_{B}, \xi])_{\lambda}(e_{1}, e_{2})&=(\mathfrak{m}_{B})_{\lambda}
(\xi(e_{1}), e_{2}))+(\mathfrak{m}_{B})_{\lambda}(e_{1}, \xi(e_{2}))
-\xi((\mathfrak{m}_{B})_{\lambda}(e_{1}, e_{2}))\\
&=-\xi(b_{1}\oo{\lambda}b_{2}),
\end{align*}
we obtain
$$
\bar{d}_{1}(\xi)_{\lambda}(e_{1}, e_{2})=\xi(b_{1})\oo{\lambda}a_{2}
+a_{1}\oo{\lambda}\xi(b_{2})-\xi(b_{1}\oo{\lambda}b_{2}),
$$
and $[\xi, \bar{d}_{1}(\xi)]=[\xi,\; \mathfrak{m}_{A}\circ(\xi, -)
+\mathfrak{m}_{A}\circ(-, \xi)-\xi\circ\mathfrak{m}_{B}]$. By direct calculation, we have
\begin{align*}
&\quad([\xi,\; \mathfrak{m}_{A}\circ(\xi, -)])_{\lambda}(e_{1}, e_{2})\\
&=\xi((\mathfrak{m}_{A})_{\lambda}\circ(\xi, -)(e_{1}, e_{2}))
-(\mathfrak{m}_{A})_{\lambda}\circ(\xi, -)(\xi(e_{1}), e_{2})
-(\mathfrak{m}_{A})_{\lambda}\circ(\xi, -)(e_{1}, \xi(e_{2}))\\
&=-\xi(b_{1})\oo{\lambda}\xi(b_{2}),\\
&\quad([\xi,\; \mathfrak{m}_{A}\circ(-, \xi)])_{\lambda}(e_{1}, e_{2})\\
&=\xi((\mathfrak{m}_{A})_{\lambda}\circ(-, \xi)(e_{1}, e_{2}))
-(\mathfrak{m}_{A})_{\lambda}\circ(-, \xi)(\xi(e_{1}), e_{2})
-(\mathfrak{m}_{A})_{\lambda}\circ(-, \xi)(e_{1}, \xi(e_{2}))\\
&=-\xi(b_{1})\oo{\lambda}\xi(b_{2}),
\end{align*}
and $([\xi, \xi\circ\mathfrak{m}_{B}])_{\lambda}(e_{1}, e_{2})=0$. That is to say,
$$
([\xi, \bar{d}_{1}(\xi)])_{\lambda}(e_{1}, e_{2})=-2\xi(b_{1})\oo{\lambda}\xi(b_{2}).
$$
Note that $(\ad_{\xi})^{n}=0$ for all $n\geq2$, since $\xi$ takes values in $A$, and
$$
\frac{e^{\ad_{\xi}}-1}{\ad_{\xi}}\bar{d}_{1}(\xi)
=\sum_{n\geq0}\frac{1}{(n+1)!}(\ad_{\xi})^{n}\bar{d}_{1}(\xi),
$$
we get
\begin{align*}
\frac{e^{\ad_{\xi}}-1}{\ad_{\xi}}\bar{d}_{1}(\xi)_{\lambda}(e_{1}, e_{2})
=&\bar{d}_{1}(\xi)_{\lambda}(e_{1}, e_{2})
+(\mbox{$\frac{1}{2}$}[\xi, \bar{d}_{1}(\xi)])_{\lambda}(e_{1}, e_{2})\\
=&\xi(b_{1})\oo{\lambda}a_{2}+a_{1}\oo{\lambda}\xi(b_{2})-\xi(b_{1}\oo{\lambda}b_{2})
-\xi(b_{1})\oo{\lambda}\xi(b_{2}).
\end{align*}
Therefore an element $\tilde{\mathfrak{c}}:=\tilde{\chi}+\tilde{\blacktriangleright}
+\tilde{\blacktriangleleft}\in{\rm MC}(L)$ is equivalent to $\mathfrak{c}$ if
\begin{align*}
\tilde{\mathfrak{c}}_{\lambda}(e_{1}, e_{2})
&=(\tilde{\chi}+\tilde{\blacktriangleright}
+\tilde{\blacktriangleleft})_{\lambda}(e_{1}, e_{2})\\
&=\mathfrak{c}_{\lambda}(e_{1}, e_{2})-b_{1}\obl{\lambda}\xi(b_{2})-\xi(b_{1})\obr{\lambda}b_{2}
-\xi(b_{1})\oo{\lambda}a_{2}\\
&\quad-a_{1}\oo{\lambda}\xi(b_{2})+\xi(b_{1}\oo{\lambda}b_{2})
+\xi(b_{1})\oo{\lambda}\xi(b_{2})\\
&=\Big(b_{1}\obl{\lambda}a_{2}-\xi(b_{1})\oo{\lambda}a_{2}\Big)
+\Big(a_{1}\obr{\lambda}b_{2}-a_{1}\oo{\lambda}\xi(b_{2})\Big)\\
&\quad+\Big(\chi_{\lambda}(b_{1}, b_{2})-b_{1}\obl{\lambda}\xi(b_{2})
-\xi(b_{1})\obr{\lambda}b_{2}+\xi(b_{1}\oo{\lambda}b_{2})+\xi(b_{1})\oo{\lambda}\xi(b_{2})\Big).
\end{align*}
Thus we get that two non-abelian 2-cocycles $(\blacktriangleright, \blacktriangleleft, \chi)$
and $(\tilde{\blacktriangleright}, \tilde{\blacktriangleleft}, \tilde{\chi})$ are equivalent,
i.e., equations (\ref{coh6})-(\ref{coh8}) are satisfied, if and only if the equation above is
satisfied, if and only if $\mathfrak{c}$ and $\tilde{\mathfrak{c}}$ are gauge equivalent.
\end{proof}

Thus, for any associative conformal algebras $A$ and $B$, we have the following
1-1 correspondences
\begin{align*}
\Ext_{nab}(B, A)\quad&\longleftrightarrow\quad \HH^{2}_{nab}(B, A)\quad
\longleftrightarrow\quad \mathcal{MC}(\mathcal{L})\\
[\mathcal{E}_{(\blacktriangleright, \blacktriangleleft, \chi)}]\quad &\;\;\leftrightarrow
\;\;\quad[(\blacktriangleright, \blacktriangleleft, \chi)]\quad \;\;\leftrightarrow \;\;\quad
[\chi+\blacktriangleright+\blacktriangleleft],
\end{align*}
where $[\mathcal{E}_{(\blacktriangleright, \blacktriangleleft, \chi)}]$ means
the equivalent class of $\mathcal{E}_{(\blacktriangleright, \blacktriangleleft, \chi)}$.

Finally, we consider the abelian extensions of associative conformal algebras.
Let $B$ be an associative conformal algebra and $A$ be an associative conformal algebra
with trivial multiplication. Then a non-abelian extension
$\mathcal{E}: \xymatrix@C=0.5cm{0\ar[r]& A\ar[r]&E\ar[r]&B\ar[r]&0}$ is
just an abelian extension of $B$ by $A$. Denote the corresponding non-abelian 2-cocycle
by $(\blacktriangleright, \blacktriangleleft, \chi)$. Then the differential graded Lie algebra
$(\mathcal{L},\; [-,-],\; \bar{d}_{\ast})$ degenerates to $(\mathcal{C}^{\bullet+1}
(B, A),\; [-,-]_{G},\, d_{\ast})$, where the bracket $[-, -]_{G}$ is trivial and
$d_{\ast}$ is given in Section \ref{sec:prel}. In this case, the Maurer-Cartan element of
$(\mathcal{C}^{\bullet+1}(B, A),\; [-,-]_{G},\, d_{\ast})$ corresponding to
extension $\mathcal{E}$ is just $\chi\in\mathcal{Z}^{2}(B, A)$, and two extensions
$\mathcal{E}$ and $\bar{\mathcal{E}}$ are equivalent if and only if
$\chi-\bar{\chi}\in\mathcal{B}^{2}(B, A)$ by equation (\ref{coh8}).
Thus, in this case, $\mathcal{MC}(\mathcal{L})$ is exactly $\HH^{2}(B, A)$,
and so that Theorem \ref{thm1} and Theorem \ref{thm2} are the same.

%%%%%%%%%%%%%%%%%%%%%%%%%%%%%%%%%%%%%%%%%%%%%%%%%%%%%%%%%%%%%%%%%%%%%%%%%%%%%%%%%
%    section  6  The inducibility problem and the fundamental exact sequence of Wells
%%%%%%%%%%%%%%%%%%%%%%%%%%%%%%%%%%%%%%%%%%%%%%%%%%%%%%%%%%%%%%%%%%%%%%%%%%%%%%%%%%%%%%

\section{The inducibility problem of automorphisms}\label{sec:ind-aut}

In this section, we study the inducibility of a pair of automorphisms about a non-abelian
extension of associative conformal algebras, and give the fundamental sequence of Wells
in the context of associative conformal algebras.

Let $A$ and $B$ be two associative conformal algebras, $\mathcal{E}: \xymatrix@C=0.5cm
{0\ar[r]&A\ar[r]^{\alpha}&E\ar[r]^{\beta}&B\ar[r]&0}$ be a non-abelian extension of $B$ by $A$
with a section $\gamma$.
Denote the corresponding non-abelian 2-cocycle of $\mathcal{E}$ by
$(\blacktriangleright, \blacktriangleleft, \chi)$. From the conclusion of Theorem \ref{thm1},
one can often identify $E$ with $A\oplus B$ under the multiplication defined
in Lemma \ref{lem: equiv}, identify $\alpha$ and $\beta$ with an injection and a projection
respectively. Denote
$$
\Aut_{A}(E):=\{f\in\Aut(E)\mid f(A)=A\}.
$$
Then, $f|_{A}\in\Aut(A)$ if $f\in\Aut_{A}(E)$. For any $f\in\Aut_{A}(E)$, we define a
$\cb[\partial]$-module homomorphism
$$
\hat{f}:=\beta\circ f\circ\gamma:\; B\longrightarrow B.
$$
Then one can check that $\hat{f}$ is independent of the choice of the section $\gamma$.
Moreover, since $\beta$ is regarded as the projection onto $B$ and $f$ preserves $B$,
we get $\hat{f}$ is a bijection. For any $b_{1}, b_{2}\in B$, we have
\begin{align*}
\hat{f}(b_{1}\oo{\lambda}b_{2})&=\beta\circ f\circ\gamma(b_{1}\oo{\lambda}b_{2})\\
&=\beta\circ f\Big(\gamma(b_{1})\oo{\lambda}\gamma(b_{2})-\chi_{\lambda}(b_{1}, b_{2})\Big)\\
&=\beta\circ f(\gamma(b_{1})\oo{\lambda}\gamma(b_{2}))\\
&=\hat{f}(b_{1})\oo{\lambda}\hat{f}(b_{2}).
\end{align*}
This means $\hat{f}\in\Aut(B)$. Note that $\widehat{g\circ f}=\beta\circ g\circ f\circ\gamma
=\beta\circ g\circ\gamma\circ\beta\circ f\circ\gamma=\hat{g}\circ\hat{f}$
for $f, g\in\Aut_{A}(E)$, we obtain a group homomorphism
$$
\kappa:\; \Aut_{A}(E)\rightarrow\Aut(A)\times\Aut(B),\qquad\qquad
f\mapsto(f|_{A}, \hat{f}).
$$

\begin{defi}\label{def:induc}
Let $A$ and $B$ be two associative conformal algebras, $\mathcal{E}: \xymatrix@C=0.5cm
{0\ar[r]&A\ar[r]^{\alpha}&E\ar[r]^{\beta}&B\ar[r]&0}$ be a non-abelian extension of $B$ by $A$,
$(\blacktriangleright, \blacktriangleleft, \chi)$ be the corresponding non-abelian 2-cocycle.
A pair $(g, h)\in\Aut(A)\times\Aut(B)$ of automorphisms is said to be {\rm inducible}
if there exists a map $f\in\Aut_{A}(E)$ such that $(g, h)=(f|_{A}, \hat{f})$.
\end{defi}

Given a pair $(g, h)\in\Aut(A)\times\Aut(B)$ and a non-abelian 2-cocycle
$(\blacktriangleright, \blacktriangleleft, \chi)$, we define a new triple
$(\blacktriangleright^{g,h}, \blacktriangleleft^{g,h}, \chi^{g,h})$ with
$\blacktriangleright^{g,h}: B\times A\rightarrow A[\lambda]$,
$\blacktriangleleft^{g,h}: A\times B\rightarrow A[\lambda]$ and
$\chi^{g,h}: B\times B\rightarrow A[\lambda]$ by
\begin{align*}
\bll{\lambda}{g,h}(b, a)&=g\Big(h^{-1}(b)\obl{\lambda}g^{-1}(a)\Big),\\
\brr{\lambda}{g,h}(a, b)&=g\Big(g^{-1}(a)\obr{\lambda}h^{-1}(b)\Big),\\
\chi^{g,h}_{\lambda}(b_{1}, b_{2})&=g\circ\chi_{\lambda}\Big(h^{-1}(b_{1}), h^{-1}(b_{2})\Big),
\end{align*}
for $a\in A$ and $b, b_{1}, b_{2}\in B$. Then one can check that
$(\blacktriangleright^{g,h}, \blacktriangleleft^{g,h}, \chi^{g,h})$ is also a non-abelian
2-cocycle on $B$ with values in $A$.

\begin{thm}\label{thm: induc}
Let $\mathcal{E}: \xymatrix@C=0.5cm{0\ar[r]& A\ar[r]^{\alpha}&E\ar[r]^{\beta}&B\ar[r]&0}$
be a non-abelian extension of associative conformal algebra $B$ by associative conformal
algebra $A$, $(\blacktriangleright, \blacktriangleleft, \chi)$ be the corresponding non-abelian
2-cocycle. With the above notations, a pair $(g, h)\in\Aut(A)\times\Aut(B)$
is inducible if and only if the non-abelian $2$-cocycles
$(\blacktriangleright, \blacktriangleleft, \chi)$ and
$(\blacktriangleright^{g,h}, \blacktriangleleft^{g,h}, \chi^{g,h})$ are equivalent.
\end{thm}

\begin{proof}
If the pair $(g, h)$ is inducible, there exists an element $f\in\Aut_{A}(E)$ such that
$f|_{A}=g$ and $\hat{f}=\beta\circ f\circ\gamma=h$, where $\gamma$ is a section of
$\mathcal{E}$. Note that for any $b\in B$,
$$
\beta\Big(f\circ\gamma\circ h^{-1}(b)-\gamma(b)\Big)=b-b=0,
$$
we get $f\circ\gamma\circ h^{-1}(b)-\gamma(b)\in\Ker(\beta)\cong A$.
Thus we can define a $\cb[\partial]$-module homomorphism $\omega: B\rightarrow A$,
$b\mapsto\gamma(b)-f\circ\gamma\circ h^{-1}(b)$. Therefore,
for any $a\in A$ and $b, b_{1}, b_{2}\in B$, we have
\begin{align*}
(\chi_{\lambda}-\chi_{\lambda}^{g,h})(b_{1}, b_{2})
&=\chi_{\lambda}(b_{1}, b_{2})-g\circ\chi_{\lambda}\Big(h^{-1}(b_{1}), h^{-1}(b_{2})\Big)\\
&=-f\Big(\gamma\circ h^{-1}(b_{1})\oo{\lambda}\gamma\circ h^{-1}(b_{2})-
\gamma(h^{-1}(b_{1})\oo{\lambda}h^{-1}(b_{2}))\Big)\\
&\qquad+\gamma(b_{1})\oo{\lambda}\gamma(b_{2})-\gamma(b_{1}\oo{\lambda}b_{2})\\
&=f\circ\gamma\circ h^{-1}(b_{1})\oo{\lambda}(\gamma(b_{2})
-f\circ\gamma\circ h^{-1}(b_{2}))\\
&\qquad+(\gamma(b_{1})-f\circ\gamma\circ h^{-1}(b_{1}))\oo{\lambda}
f\circ\gamma\circ h^{-1}(b_{2})\\
&\qquad-(\gamma(b_{1})-f\circ\gamma\circ h^{-1}(b_{1}))\oo{\lambda}
(\gamma(b_{2})-f\circ\gamma\circ h^{-1}(b_{2}))\\
&\qquad-(\gamma-f\circ\gamma\circ h^{-1})(b_{1}\oo{\lambda}b_{2})\\
&=b_{1}\obl{\lambda}\omega(b_{2})+\omega(b_{1})\obr{\lambda}b_{2}
-\omega(b_{1}\oo{\lambda}b_{2})-\omega(b_{1})\oo{\lambda}\omega(b_{2}),
\end{align*}
and
\begin{align*}
(\obl{\lambda}-\bll{\lambda}{g,h})(b, a)
&=b\obl{\lambda}a-g\Big(h^{-1}(b)\obl{\lambda}g^{-1}(a)\Big)\\
&=\gamma(b)\oo{\lambda}a-f\Big(\gamma\circ h^{-1}(b)\oo{\lambda}g^{-1}(a)\Big)\\
&=\Big(\gamma(b)-f\circ\gamma\circ h^{-1}(b)\Big)\oo{\lambda}a\\
&=\omega(b)\oo{\lambda}a.
\end{align*}
Similarly, $(\obr{\lambda}-\brr{\lambda}{g,h})(a, b)=a\oo{\lambda}\omega(b)$.
Thus, it follows from Definition \ref{def:coho} that $(\blacktriangleright,
\blacktriangleleft, \chi)$ and $(\blacktriangleright^{g,h}, \blacktriangleleft^{g,h},
\chi^{g,h})$ are equivalent.

On the other hand, if the non-abelian $2$-cocycles
$(\blacktriangleright, \blacktriangleleft, \chi)$ and
$(\blacktriangleright^{g,h}, \blacktriangleleft^{g,h}, \chi^{g,h})$ are equivalent,
there exists a $\cb[\partial]$-module homomorphism $\omega: B\rightarrow A$
satisfying the equations (\ref{coh6})-(\ref{coh8}), and there is a commutative diagram:
$$
\xymatrix{0\ar[r]& A\ar[r]^{\alpha\qquad}\ar@{=}[d]&A\oplus_{(\blacktriangleright,
\blacktriangleleft, \chi)}B\ar[r]^{\qquad\beta}\ar[d]^{\theta}
&B\ar[r]\ar@{=}[d]&0\;\\ 0\;\;\ar[r]& A\ar[r]^{\bar{\alpha}\qquad}
&A\oplus_{(\blacktriangleright^{g,h}, \blacktriangleleft^{g,h}, \chi^{g,h})}
B\ar[r]^{\quad\qquad\bar{\beta}}&B\ar[r]&0,}
$$
where $\theta(a, b)=(a-\omega(b),\; b)$. Since $E\cong A\oplus_{(\blacktriangleright,
\blacktriangleleft, \chi)}B$ as associative conformal algebras, we denote each element
in $E$ by $(a, b)$, and define a $\cb[\partial]$-module homomorphism
$f: E\rightarrow E$ by $f(a, b)=\theta(g(a), h(b))=(g(a)+\omega(h(b)),\; h(b))$, for any
$a\in A$ and $b\in B$. Then $f$ is a bijection since $\theta$ is an isomorphism.
Moreover, for any $(a_{1}, b_{1})$, $(a_{2}, b_{2})\in E$, we have
\begin{align*}
&f(a_{1}, b_{1})\oo{\lambda}f(a_{2}, b_{2})\\
=&\; \theta\Big((g(a_{1}), h(b_{1}))\oo{\lambda}(g(a_{2}), h(b_{2}))\Big)\\
=&\; \theta\Big(g(a_{1})\oo{\lambda}g(a_{2})+h(b_{1})\bll{\lambda}{g,h}g(a_{2})
+g(a_{1})\brr{\lambda}{g,h}h(b_{2})+\chi_{\lambda}^{g,h}(h(b_{1}), h(b_{2})),\;
h(b_{1})\oo{\lambda}h(b_{2})\Big)\\
=&\; f\Big(a_{1}\oo{\lambda}a_{2}+b_{1}\obl{\lambda}a_{2}+a_{1}\obr{\lambda}b_{1}
+\chi_{\lambda}(b_{1}, b_{2}),\; b_{1}\oo{\lambda}b_{2}\Big)\\
=&\; f\Big((a_{1}, b_{1})\oo{\lambda}(a_{2}, b_{2})\Big).
\end{align*}
This means $f\in\Aut(E)$. Clearly, $f|_{A}=g$ and $\hat{f}=\beta\circ f\circ\gamma=h$.
Hence the pair $(g, h)$ is inducible.
\end{proof}

From the equations (\ref{coh6})-(\ref{coh8}), we can directly get the following corollary.

\begin{cor}\label{cor:describe}
Let $\mathcal{E}: \xymatrix@C=0.5cm{0\ar[r]& A\ar[r]&E\ar[r]&B\ar[r]&0}$
be a non-abelian extension of associative conformal algebra $B$ by associative conformal
algebra $A$, $(\blacktriangleright, \blacktriangleleft, \chi)$ be the corresponding non-abelian
2-cocycle. A pair $(g, h)\in\Aut(A)\times\Aut(B)$ is inducible if and
only if there exists a $\cb[\partial]$-module homomorphism $\omega: B\rightarrow A$
such that for any $a\in A$ and $b, b_{1}, b_{2}\in B$,
\begin{align*}
g(b\obl{\lambda}a)&=h(b)\obl{\lambda}g(a)-\omega\circ h(b)\oo{\lambda}g(a),\\
g(a\obr{\lambda}b)&=g(a)\obr{\lambda}h(b)-g(a)\oo{\lambda}\omega\circ h(b),\\
g\circ\chi_{\lambda}(b_{1}, b_{2})&=\chi_{\lambda}(h(b_{1}), h(b_{2}))
-h(b_{1})\obl{\lambda}\omega\circ h(b_{2})+\omega(h(b_{1})\oo{\lambda}h(b_{2}))\\
&\qquad-\omega\circ h(b_{1})\obr{\lambda}h(b_{2})
+\omega\circ h(b_{1})\oo{\lambda}\omega\circ h(b_{2}).
\end{align*}
\end{cor}

Next, we will interpret the above theorem in terms of the Wells map in the context of
associative conformal algebras. Let $\mathcal{E}: \xymatrix@C=0.5cm{0\ar[r]&
A\ar[r]&E\ar[r]&B\ar[r]&0}$ be a non-abelian extension of $B$ by $A$, $(\blacktriangleright,
\blacktriangleleft, \chi)$ be the corresponding non-abelian 2-cocycle.
We define a map $\mathcal{W}:\; \Aut(A)\times\Aut(B)\rightarrow \HH^2_{nab}(B, A)$ by
$$
\mathcal{W}(g, h)=[(\blacktriangleright^{g,h}, \blacktriangleleft^{g,h}, \chi^{g,h})
-(\blacktriangleright, \blacktriangleleft, \chi)].
$$
Then one can check that $\mathcal{W}$ is a group homomorphism, which is called {\it Wells map}.
Therefore, we can directly obtain the following corollary by Theorem \ref{thm: induc}.

\begin{cor}\label{cor:wells}
Let $\mathcal{E}: \xymatrix@C=0.5cm{0\ar[r]& A\ar[r]&E\ar[r]&B\ar[r]&0}$
be a non-abelian extension of associative conformal algebra $B$ by associative conformal
algebra $A$. A pair $(g, h)\in\Aut(A)\times\Aut(B)$ of automorphisms is inducible if and
only if $\mathcal{W}(g, h)=0$.
\end{cor}

Let $\mathcal{E}: \xymatrix@C=0.5cm{0\ar[r]&A\ar[r]&E\ar[r]&B\ar[r]&0}$ be a non-abelian
extension of $B$ by $A$. Next, we define
$$
\Aut_{A}^{\id}(E):=\{f\in\Aut_{A}(E)\mid\kappa(f)=(\id_{A}, \id_{B})\}.
$$
Then $\Aut_{A}^{\id}(E)$ is a subgroup of $\Aut_{A}(E)$ and there exists an exact
sequence as following.

\begin{cor}\label{cor: wells}(Fundamental sequence of Wells)
With the above notations, there is an exact sequence
$$
\xymatrix{1\ar[r]&\Aut_{A}^{\id}(E)\ar[r]^{\iota}&\Aut_{A}(E)\ar[r]^{\;\kappa\qquad}
&\Aut(A)\times\Aut(B)\ar[r]^{\quad\mathcal{W}\;}&\HH^{2}_{nab}(B, A)},
$$
where $\iota$ is the inclusion map.
\end{cor}

\begin{proof}
We need show this sequence is exact at $\Aut_{A}(E)$ and $\Aut(A)\times\Aut(B)$.
First, note that an element $f\in\Aut_{A}(E)$ such that $f\in\Ker(\kappa)$
if and only if $\kappa(f)=(\id_{A}, \id_{B})$, i.e., $f\in\Aut_{A}^{\id}(E)$,
we obtain the sequence is exact at $\Aut_{A}(E)$. Second, by Corollary \ref{cor:wells},
we get a pair $(g, h)\in\Aut(A)\times\Aut(B)$ satisfying $(g, h)\in\Img(\kappa)$, i.e.,
$(g, h)$ is inducible, if and only if $\mathcal{W}(g, h)=0$. This means the sequence
is exact at $\Aut(A)\times\Aut(B)$.
\end{proof}

More generally, if define
\begin{align*}
\Aut_{A}^{A}(E):&=\{f\in\Aut_{A}(E)\mid f|_{A}=\id_{A}\},\\
\Aut_{A}^{B}(E):&=\{f\in\Aut_{A}(E)\mid\hat{f}=\id_{B}\},
\end{align*}
we obtain two morphisms of groups $\kappa_{B}: \Aut_{A}^{A}(E)\rightarrow\Aut(B)$,
$f\mapsto\hat{f}$, and $\kappa_{A}: \Aut_{A}^{B}(E)\rightarrow\Aut(A)$,
$f\mapsto f|_{A}$. Define maps $\mathcal{W}_{A}: \Aut(A)\rightarrow\HH^{2}_{nab}(B, A)$
and $\mathcal{W}_{B}: \Aut(B)\rightarrow\HH^{2}_{nab}(B, A)$ by
\begin{align*}
\mathcal{W}_{A}(g):&=[(\blacktriangleright^{g,\id_{B}}, \blacktriangleleft^{g,\id_{B}},
\chi^{g,\id_{B}})-(\blacktriangleright, \blacktriangleleft, \chi)],\\
\mathcal{W}_{B}(h):&=[(\blacktriangleright^{\id_{A},h}, \blacktriangleleft^{\id_{A},h},
\chi^{\id_{A},h})-(\blacktriangleright, \blacktriangleleft, \chi)].
\end{align*}
Then similar to Corollary \ref{cor: wells}, we have two exact sequences of groups.

\begin{pro}\label{pro: g-wells}
Let $\mathcal{E}: \xymatrix@C=0.5cm{0\ar[r]&A\ar[r]&E\ar[r]&B\ar[r]&0}$ be a non-abelian
extension of associative conformal algebra $B$ by associative conformal algebra $A$.
With the above notations, there are two exact sequences
\begin{align*}
&\xymatrix{1\ar[r]&\Aut_{A}^{\id}(E)\ar[r]^{\iota_{A}}&\Aut_{A}^{B}(E)\ar[r]^{\;\kappa_{A}}
&\Aut(A)\ar[r]^{\mathcal{W}_{A}\quad}&\HH^{2}_{nab}(B, A)},\\
&\xymatrix{1\ar[r]&\Aut_{A}^{\id}(E)\ar[r]^{\iota_{B}}&\Aut_{A}^{A}(E)\ar[r]^{\;\kappa_{B}}
&\Aut(B)\ar[r]^{\mathcal{W}_{B}\quad}&\HH^{2}_{nab}(B, A)},
\end{align*}
where $\iota_{A}$ and $\iota_{B}$ are inclusion maps.
\end{pro}

By using the exact sequences in this proposition, we can give some results of
the lifting problem of automorphism in a non-abelian extension.

\begin{cor}\label{cor:lift}
Let $\mathcal{E}: \xymatrix@C=0.5cm{0\ar[r]& A\ar[r]&E\ar[r]&B\ar[r]&0}$
be a non-abelian extension of associative conformal algebra $B$ by associative conformal
algebra $A$.

$(i)$ An automorphism $g\in\Aut(A)$ can be extended to an automorphism of $E$ inducing the
identity on $B$ if and only if $\mathcal{W}_{A}(g)=0$.

$(ii)$ An automorphism $h\in\Aut(B)$ can be lifted to an automorphism of $E$ fixing
$B$ pointwise if and only if $\mathcal{W}_{B}(h)=0$.
\end{cor}

At the end of this section, let's go back to the abelian extension.
Let $B$ be an associative conformal algebra, $A$ be a bimodule over $B$, and
$\mathcal{E}: \xymatrix@C=0.5cm{0\ar[r]& A\ar[r]^{\alpha}&E\ar[r]^{\beta}&B\ar[r]&0}$
be an abelian extension of $B$ by $A$, where the $B$-bimodule structure on $A$
is given by $(\blacktriangleright, \blacktriangleleft)$. Denote by $\Aut(A)$ the set of
all $\cb[\partial]$-module automorphisms, and by
$\Aut_{\blacktriangleright, \blacktriangleleft}(A, B)$ the set of all pairs $(g, h)\in
\Aut(A)\times\Aut(B)$ satisfying
$$
g(b\blacktriangleright_{\lambda}a)=h(b)\blacktriangleright_{\lambda}g(a),\qquad\qquad
g(a\blacktriangleleft_{\lambda}b)=g(a)\blacktriangleleft_{\lambda}h(b),
$$
for $a\in A$ and $b\in B$. Then $\Aut_{\blacktriangleright, \blacktriangleleft}(A, B)$
is a subgroup of $\Aut(A)\times\Aut(B)$. Let $\gamma$ be a section of $\mathcal{E}$,
and $\chi$ be the 2-cocycle corresponding to $\mathcal{E}$. For a pair $(g, h)\in
\Aut(A)\times\Aut(B)$, we define $\chi^{g,h}: B\times B\rightarrow A[\lambda]$ by
$$
\chi^{g,h}_{\lambda}(b_{1}, b_{2})=g\circ\chi_{\lambda}(h^{-1}(b_{1}), h^{-1}(b_{2})),
$$
for $b_{1}, b_{2}\in B$. If $(g, h)\in\Aut_{\blacktriangleright,\blacktriangleleft}(A, B)$,
we have
$$
d_{2}(\chi^{g,h})_{\lambda_{1}, \lambda_{2}}(b_{1}, b_{2}, b_{3})
=g\circ d_{2}(\chi)_{\lambda_{1}, \lambda_{2}}\Big(h^{-1}(b_{1}), h^{-1}(b_{2}),
h^{-1}(b_{3})\Big)=0,
$$
for any $b_{1}, b_{2}, b_{3}\in B$, since $d_{2}(\chi)=0$. That is, $\chi^{g,h}
\in\mathcal{Z}^{2}(B, A)$. From Theorem \ref{thm: induc} and Corollary \ref{cor:describe},
we obtain following theorem.

\begin{thm}\label{thm: abe-induc}
Let $\mathcal{E}: \xymatrix@C=0.5cm{0\ar[r]&A\ar[r]&E\ar[r]&B\ar[r]&0}$ be an abelian
extension of associative conformal algebra $B$ by bimodule $A$.
Then a pair $(g, h)\in\Aut_{\blacktriangleright,\blacktriangleleft}(A, B)$ is inducible
if and only if the 2-cocycles $\chi$ and $\chi^{g,h}$ are equivalent, if and
only if there exists a $\cb[\partial]$-module homomorphism $\omega: B\rightarrow A$
such that for any $b_{1}, b_{2}\in B$,
\begin{align*}
g\circ\chi_{\lambda}(b_{1}, b_{2})&=\chi_{\lambda}(h(b_{1}), h(b_{2}))
-h(b_{1})\obl{\lambda}\omega\circ h(b_{2})\\
&\qquad+\omega(h(b_{1})\oo{\lambda}h(b_{2}))-\omega\circ h(b_{1})\obr{\lambda}h(b_{2}).
\end{align*}
\end{thm}

We next define the Wells map for the abelian extensions by
$$
\bar{\mathcal{W}}:\; \Aut_{\blacktriangleright,\blacktriangleleft}(A, B)\rightarrow
\HH^{2}(B, A)\qquad\qquad (g, h)\mapsto[\chi^{g, h}-\chi].
$$
Let $\gamma$ be a section of $\mathcal{E}$. For any $f\in\Aut_{A}(E)$, we define
$\hat{f}:=\beta\circ f\circ\gamma$. Then we obtain a group homomorphism
$$
\bar{\kappa}:\; \Aut_{A}(E)\rightarrow\Aut_{\blacktriangleright,\blacktriangleleft}(A, B),
\qquad\qquad f\mapsto(f|_{A}, \hat{f}).
$$
The fundamental sequence of Wells can be obtained as the following form.

\begin{cor}\label{cor:abe-wells}
Let $\mathcal{E}: \xymatrix@C=0.5cm{0\ar[r]& A\ar[r]&E\ar[r]&B\ar[r]&0}$
be an abelian extension of associative conformal algebra $B$ by bimodule $A$.
With the above notations, there is an exact sequence
$$
\xymatrix{1\ar[r]&\Aut_{A}^{\id}(E)\ar[r]^{\iota}&\Aut_{A}(E)\ar[r]^{\bar{\kappa}\quad}
&\Aut_{\blacktriangleright,\blacktriangleleft}(A, B)\ar[r]^{\bar{\mathcal{W}}} &\HH^{2}(B, A)},
$$
where $\iota$ is the inclusion map.
\end{cor}

In particular, if the abelian extension $\mathcal{E}: \xymatrix@C=0.5cm{0\ar[r]&
A\ar[r]^{\alpha}&E\ar[r]^{\beta}&B\ar[r]&0}$ is split in the category of associative
conformal algebras, then the associative conformal algebra $E$ is isomorphic to the
semidirect product $A\rtimes B$, where the conformal multiplication on $A\oplus B$ is given by
$$
(a_{1}, b_{1})\oo{\lambda}(a_{2}, b_{2})=(a_{1}\obr{\lambda}b_{2}
+b_{1}\obl{\lambda}a_{2},\; b_{1}\oo{\lambda}b_{2}),
$$
for $(a_{1}, b_{1})$, $(a_{2}, b_{2})\in A\oplus B$, and the section $\gamma: B\rightarrow E$,
is given by $\gamma(b)=(0, b)$ for $b\in B$. In this case the corresponding 2-cocycle
$\chi$ of $\mathcal{E}$ is zero, and the Wells map vanishes identically.
Thus we have an exact sequence
$$
\xymatrix{1\ar[r]&\Aut_{A}^{\id}(E)\ar[r]^{\iota}&\Aut_{A}(E)\ar[r]^{\bar{\kappa}\quad}
&\Aut_{\blacktriangleright,\blacktriangleleft}(A, B)\ar[r]&1}.
$$
We define map $\varrho: \Aut_{\blacktriangleright,\blacktriangleleft}(A, B)\rightarrow
\Aut_{A}(E)$ by
$$
\varrho(g,h)(a, b)=(g(a), h(b)),
$$
for any $(g, h)\in\Aut_{\blacktriangleright,\blacktriangleleft}(A, B)$ and $(a, b)\in E$.
Then $\bar{\kappa}\circ\varrho(g,h)=(\varrho(g,h)|_{A},\; \widehat{\varrho(g,h)})=(g, h)$,
since $\beta\circ\varrho(g,h)\circ\gamma(b)=h(b)$ for any $b\in B$.
That is to say, the exact sequence above is split in the category of groups.
Thus we have the following corollary.

\begin{cor}\label{cor:split}
Let $\mathcal{E}: \xymatrix@C=0.5cm{0\ar[r]& A\ar[r]&E\ar[r]&B\ar[r]&0}$ be a split abelian
extension of associative conformal algebra $B$ by bimodule $A$. Then
$$
\Aut_{A}(E)\cong \Aut_{A}^{\id}(E)\times\Aut_{\blacktriangleright,\blacktriangleleft}(A, B)
$$
as groups.
\end{cor}

%%%%%%%%%%%%%%%%%%%%%%%%%%%%%%%%%%%%%%%%%%%%%%%%%%%%%%%%%%%%%%%%%%%%%%%%%%%%%%%%%
%    section  7  The extensibility problem
%%%%%%%%%%%%%%%%%%%%%%%%%%%%%%%%%%%%%%%%%%%%%%%%%%%%%%%%%%%%%%%%%%%%%%%%%%%%%%%%%%%%%%

\section{The extensibility problem of derivations}\label{sec:ind-der}

For an abelian extension of associative conformal algebras, we construct a
Lie algebra such that the 2-th cohomology group is a representation over it.
Using this representation, we study the extensibility of a pair of derivations
about an abelian extension of associative conformal algebras, and give an exact sequence
of Wells type.

Let $B$ be an associative conformal algebra. Recall that a {\it derivation} of $B$ is
a $\cb[\partial]$-module homomorphism $\hd: B\rightarrow B$ such that
$\hd(b_{1}\oo{\lambda}b_{2})=\hd(b_{1})\oo{\lambda}b_{2}+b_{1}\oo{\lambda}\hd(b_{2})$,
for any $b_{1}, b_{1}\in B$, where the first $\hd$ is extended canonically to a
$\cb[\partial]$-module homomorphism from $B[\lambda]$ to $B[\lambda]$, i.e.,
$\hd(\sum b_{i}\lambda^{i})=\sum \hd(b_{i})\lambda^{i}$. We denote the set of all
derivations of $B$ by $\Der(B)$. Then $\Der(B)$ is a Lie algebra over $\cb$ with respect to
the commutator. Let $A$ be a (conformal) $B$-bimodule, and
$$
\mathcal{E}:\qquad \xymatrix{0\ar[r]&A\ar[r]^{\alpha}&E\ar[r]^{\beta}&B\ar[r]&0}
$$
be an abelian extension of $B$ by $A$. Then there is an associative
conformal algebra structure on $\cb[\partial]$-module $A\oplus B$, which is given by
$$
(a_{1}, b_{1})\oo{\lambda}(a_{2}, b_{2})=(a_{1}\oor{\lambda}b_{2}
+b_{1}\ool{\lambda}a_{2},\; b_{1}\oo{\lambda}b_{2}),
$$
for any $a_{1}, a_{2}\in A$ and $b_{1}, b_{2}\in B$, where $(\triangleright, \triangleleft)$
is the $B$-bimodule action on $A$. Here the module action of $B$ on $A$ is also given by the
multiplication of $E$. Indeed, let $\gamma: B\rightarrow E$ be a $\cb[\partial]$-module
homomorphism such that $\beta\circ\gamma=\id_{B}$. Then $a_{1}\oor{\lambda}b_{2}
=a_{1}\oo{\lambda}\gamma(b_{2})$, $b_{1}\ool{\lambda}a_{2}=\gamma(b_{1})\oo{\lambda}a_{2}$,
and this action does not depend on the choice of $\gamma$. We denote this associative conformal
algebra by $A\rtimes B$. We view $A$ as a trivial associative conformal algebra.
Then $\Der(A)$ just the set of all $\cb[\partial]$-module endomorphisms of $A$.
Then by direct calculation, we have the following lemma.

\begin{lem} \label{lem: der}
Let $(\hd_{A}, \hd_{B})\in\Der(A)\times\Der(B)$. Then $(\hd_{A}, \hd_{B})\in\Der(A\rtimes B)$
if and only if
\begin{align}
\hd_{A}(b\ool{\lambda}a)&=b\ool{\lambda}\hd_{A}(a)+\hd_{B}(b)\ool{\lambda}a,\label{der1}\\
\hd_{A}(a\oor{\lambda}b)&=a\oor{\lambda}\hd_{B}(b)+\hd_{A}(a)\oor{\lambda}b,\label{der2}
\end{align}
for any $a\in A$ and $b\in B$.
\end{lem}

Denote by
$$
\g(A, B):=\Big\{(\hd_{A}, \hd_{B})\in\Der(A)\times\Der(B)\mid(\hd_{A}, \hd_{B})
\mbox{ satisfies equations (\ref{der1}) and (\ref{der2})} \Big\}.
$$
Then one can check that $\g(A, B)$ is a Lie subalgebra of $\Der(A\rtimes B)$.
Now we define an action of $\g(A, B)$ on the space of $2$-cochains $\mathcal{C}^{2}(B, A)$ by
$$
\Theta(\hd_{A}, \hd_{B})(f)_{\lambda}(b_{1}, b_{2})=\hd_{A}(f_{\lambda}(b_{1}, b_{2}))
-f_{\lambda}(\hd_{B}(b_{1}), b_{2})-f_{\lambda}(b_{1}, \hd_{B}(b_{2})),
$$
for any $(\hd_{A}, \hd_{B})\in\g(A, B)$, $f\in\mathcal{C}^{2}(B, A)$ and
$b_{1}, b_{2}\in B$. Then $d_{2}(\Theta(\hd_{A}, \hd_{B})(f))=0$
for any $f\in\mathcal{Z}^{2}(B, A)$, and $\Theta(\hd_{A}, \hd_{B})(d_{1}(f))
=d_{1}(\hd_{A}\circ f-f\circ\hd_{B})$ for any $f\in\mathcal{C}^{1}(B, A)$.
Thus we get a linear map
$$
\Theta: \g(A, B)\longrightarrow\gl(\HH^{2}(B, A)),\qquad\qquad
\Theta(\hd_{A}, \hd_{B})([f])=[\Theta(\hd_{A}, \hd_{B})(f)],
$$
for any $(\hd_{A}, \hd_{B})\in\g(A, B)$ and $[f]\in\HH^{2}(B, A)$.
Moreover, one can show that $\Theta$ is a homomorphism of Lie algebras as following.

\begin{pro}\label{pro:rep-lie}
Then map $\Theta$ gives a representation of Lie algebra $\g(A, B)$ on $\HH^{2}(B, A)$.
\end{pro}

\begin{proof}
For any $(\hd_{A}, \hd_{B})$, $(\hd'_{A}, \hd'_{B})\in\Der(A)\times\Der(B)$ and
$[f]\in\HH^{2}(B, A)$, note that
\begin{align*}
&\; \Big(\Theta(\hd_{A}, \hd_{B})\circ\Theta(\hd'_{A}, \hd'_{B})\Big)([f])\\
=&\; \Theta(\hd_{A}, \hd_{B})\Big([\hd'_{A}\circ f-f\circ(\hd'_{B}\otimes\id_{B})-
f\circ(\id_{B}\otimes\hd'_{B})]\Big)\\
=&\; \Big[\hd_{A}\circ(\hd'_{A}\circ f-f\circ(\hd'_{B}\otimes\id_{B})-
f\circ(\id_{B}\otimes\hd'_{B}))\\
&\quad-(\hd'_{A}\circ f-f\circ(\hd'_{B}\otimes\id_{B})-
f\circ(\id_{B}\otimes\hd'_{B}))(\hd_{B}\otimes\id_{B})\\
&\quad-(\hd'_{A}\circ f-f\circ(\hd'_{B}\otimes\id_{B})-
f\circ(\id_{B}\otimes\hd'_{B}))(\id_{B}\otimes\hd_{B})\Big],
\end{align*}
we get
\begin{align*}
&\; \Big[\Theta(\hd_{A}, \hd_{B}),\; \Theta(\hd'_{A}, \hd'_{B})\Big]([f])\\
=&\; \Big(\Theta(\hd_{A}, \hd_{B})\circ\Theta(\hd'_{A}, \hd'_{B})\Big)([f])
-\Big(\Theta(\hd'_{A}, \hd'_{B})\circ\Theta(\hd_{A}, \hd_{B})\Big)([f])\\
=&\; \Theta(\hd_{A}\circ\hd'_{A}-\hd'_{A}\circ\hd_{A},\;
\hd_{B}\circ\hd'_{B}-\hd'_{B}\circ\hd_{B})([f])\\
=&\; \Theta\Big([(\hd_{A}, \hd_{B}),\; (\hd'_{A}, \hd'_{B})]\Big)([f]).
\end{align*}
Thus $\Theta$ is a homomorphism of Lie algebras.
\end{proof}

Now we consider the extensibility of a pair of derivations
about an abelian extension of associative conformal algebras.

\begin{defi}\label{def:exte}
Let $\mathcal{E}: \xymatrix@C=0.5cm{0\ar[r]& A\ar[r]^{\alpha}&E\ar[r]^{\beta}&B\ar[r]&0}$
be an abelian extension of associative conformal algebra $B$ by bimodule $A$.
A pair $(\hd_{A}, \hd_{B})\in\Der(A)\times\Der(B)$ is called an {\rm extensible
pair} if there is a derivation $\hd_{E}\in\Der(E)$ such that $\hd_{E}\circ\alpha=
\alpha\circ\hd_{A}$ and $\hd_{B}\circ\beta=\beta\circ\hd_{E}$.
\end{defi}

For the abelian extension $\mathcal{E}$, there is a $\cb[\partial]$-module homomorphism
$\gamma: B\rightarrow E$ such that $\beta\circ\gamma=\id_{B}$.
Define $\chi: B\times B\rightarrow A[\lambda]$ by
$$
\chi_{\lambda}(b_{1}, b_{2})=\gamma(b_{1})\oo{\lambda}\gamma(b_{2})
-\gamma(b_{1}\oo{\lambda} b_{2}),
$$
for any $b_{1}, b_{2}\in B$. By Section \ref{sec:Non-abel}, we get $\chi$ is a
2-cocycle in $\mathcal{Z}^{2}(B, A)$, and the cohomology class $[\chi]\in\HH^{2}(B, A)$
does not depend on the choice of $\gamma$. Thus $\Theta(\hd_{A}, \hd_{B})([\chi])$
does not depend on the choice of $\gamma$. Using $\Theta(\hd_{A}, \hd_{B})([\chi])$,
we can define a Wells map associated to $\mathcal{E}$ as follows.

\begin{defi}\label{def:well}
Let $\mathcal{E}: \xymatrix@C=0.5cm{0\ar[r]& A\ar[r]^{\alpha}&E\ar[r]^{\beta}&B\ar[r]&0}$
be an abelian extension of associative conformal algebra $B$ by bimodule $A$.
The map $\mathcal{W}: \g(A, B)\rightarrow\HH^{2}(B, A)$,
$$
\mathcal{W}(\hd_{A}, \hd_{B})=\Theta(\hd_{A}, \hd_{B})([\chi])
$$
is called the {\rm Wells map} associated to $\mathcal{E}$.
\end{defi}

If the abelian extension $\mathcal{E}$ is split, i.e., there is a homomorphism
of associative conformal algebras $\gamma: B\rightarrow E$ such that $\beta\circ\gamma
=\id_{B}$, then $\chi=0$ and $\mathcal{W}=0$. More generally, we have the following theorem.

\begin{thm}\label{thm: ext-con}
Let $\mathcal{E}: \xymatrix@C=0.5cm{0\ar[r]& A\ar[r]^{\alpha}&E\ar[r]^{\beta}&B\ar[r]&0}$
be an abelian extension of associative conformal algebra $B$ by bimodule $A$.
A pair $(\hd_{A}, \hd_{B})\in\Der(A)\times\Der(B)$ is extensible if and only if
$(\hd_{A}, \hd_{B})\in\g(A, B)$ and $\mathcal{W}(\hd_{A}, \hd_{B})=[0]$.
\end{thm}

\begin{proof}
Let $\gamma$ be a section of $\mathcal{E}$. We identify $A$ as a $\cb[\partial]$-submodule
of $E$. Then the map $\tau: A\oplus B\rightarrow E$, $(a, b)\mapsto a+\gamma(b)$
gives an isomorphism of $\cb[\partial]$-modules. For convenience, denote an element
in $E$ by $a+\gamma(b)$ for $a\in A$ and $b\in B$.

If $(\hd_{A}, \hd_{B})\in\g(A, B)$ and $\mathcal{W}(\hd_{A}, \hd_{B})=[0]$, there exists
a map $f\in\Hom_{\cb[\partial]}(B, A)$ such that $\Omega(\hd_{A}, \hd_{B})=[d_{1}(f)]$.
We define a map $\hd_{E}: E\rightarrow E$ by
$$
\hd_{E}(a+\gamma(b))=\hd_{A}(a)+f(b)+\gamma(\hd_{B}(b)),
$$
for any $a\in A$ and $b\in B$. It is easy to see that $\hd_{E}$ is a homomorphism of
$\cb[\partial]$-modules, and
\begin{align*}
\hd_{E}(\alpha(a))&=\hd_{E}(a)=\hd_{A}(a)=\alpha(\hd_{A}(a)),\\
\beta(\hd_{E}(a+\gamma(b)))&=\beta(\hd_{A}(a)+f(b)+\gamma(\hd_{B}(b)))
=\hd_{B}(b)=\hd_{B}(\beta(a+\gamma(b))),
\end{align*}
for any $a\in A$ and $b\in B$. Moreover, for any $b_{1}, b_{2}\in B$, by the definition
of $\hd_{E}$, we have
\begin{align*}
\gamma(\hd_{B}(b_{1}))\oo{\lambda}\gamma(b_{2})
&=\hd_{E}(\gamma(b_{1}))\oo{\lambda}\gamma(b_{2})-f(b_{1})\oo{\lambda}\gamma(b_{2}),\\
\gamma(b_{1})\oo{\lambda}\gamma(\hd_{B}(b_{2}))
&=\gamma(b_{1})\oo{\lambda}\hd_{E}(\gamma(b_{2}))-\gamma(b_{1})\oo{\lambda}f(b_{2}).
\end{align*}
Thus,
\begin{align*}
\hd_{E}\Big(\gamma(b_{1})\oo{\lambda}\gamma(b_{2})\Big)
&=\hd_{E}\Big(\chi_{\lambda}(b_{1}, b_{2})+\gamma(b_{1}\oo{\lambda}b_{2})\Big)\\
&=\hd_{A}(\chi_{\lambda}(b_{1}, b_{2}))+f(b_{1}\oo{\lambda}b_{2})
+\gamma\Big(\hd_{B}(b_{1})\oo{\lambda}b_{2}+b_{1}\oo{\lambda}\hd_{B}(b_{2})\Big)\\
&=\hd_{A}(\chi_{\lambda}(b_{1}, b_{2}))+f(b_{1}\oo{\lambda}b_{2})
+\gamma(\hd_{B}(b_{1}))\oo{\lambda}\gamma(b_{2})\\
&\quad -\chi_{\lambda}(\hd_{B}(b_{1}),\; b_{2})
+\gamma(b_{1})\oo{\lambda}\gamma(\hd_{B}(b_{1}))-\chi_{\lambda}(b_{1},\; \hd_{B}(b_{2}))\\
&=\hd_{E}(\gamma(b_{1}))\oo{\lambda}\gamma(b_{2})
+\gamma(b_{1})\oo{\lambda}\hd_{E}(\gamma(b_{2}))\\
&\quad +\Theta(\hd_{A}, \hd_{B})(\chi)_{\lambda}(b_{1}, b_{2})
-d_{1}(f)_{\lambda}(b_{1}, b_{2})\\
&=\hd_{E}(\gamma(b_{1}))\oo{\lambda}\gamma(b_{2})
+\gamma(b_{1})\oo{\lambda}\hd_{E}(\gamma(b_{2})).
\end{align*}
Therefore, by the equations (\ref{der1}) and (\ref{der2}), we obtain $\hd_{E}\in\Der(E)$.
Hence the pair $(\hd_{A}, \hd_{B})\in\Der(A)\times\Der(B)$ is extensible.

Conversely, if there exists a derivation $\hd_{E}\in\Der(E)$ such that $\hd_{E}\circ\alpha
=\alpha\circ\hd_{A}$ and $\beta\circ\hd_{E}=\hd_{B}\circ\beta$, then $\hd_{A}=\hd_{E}|_{A}$.
Since $\beta\Big(\hd_{E}(\gamma(b))-\gamma(\hd_{B}(b))\Big)=0$ for any $b\in B$,
we get a $\cb[\partial]$-module homomorphism $\delta: B\rightarrow A$,
$$
\delta(b)=\hd_{E}(\gamma(b))-\gamma(\hd_{B}(b)),
$$
for $b\in B$. Then, for any $a\in A$ and $b\in B$,
\begin{align*}
\hd_{E}(\gamma(b)\oo{\lambda}a)-\gamma(b)\oo{\lambda}\hd_{E}(a)
&=\hd_{E}(\gamma(b))\oo{\lambda}a+\gamma(b)\oo{\lambda}\hd_{E}(a)
-\gamma(b)\oo{\lambda}\hd_{E}(a)\\
&=\Big(\gamma(\hd_{B}(b))+\delta(b)\Big)\oo{\lambda}a\\
&=\gamma(\hd_{B}(b))\oo{\lambda}a.
\end{align*}
This means the equation (\ref{der1}) holds. Similarly, one can check that
the equation (\ref{der2}) also holds. Hence $(\hd_{A}, \hd_{B})\in\g(A, B)$.
Finally, we show that $\Theta(\hd_{A}, \hd_{B})(\chi)=d_{1}(\delta)$.
For any $b_{1}, b_{2}\in B$, we have
\begin{align*}
&\; \hd_{E}\Big(\gamma(b_{1})\oo{\lambda}\gamma(b_{2})\Big)\\
=&\; \hd_{E}\Big(\gamma(b_{1}\oo{\lambda}b_{2})+\chi_{\lambda}(b_{1}, b_{2})\Big)\\
=&\; \gamma(\hd_{B}(b_{1}\oo{\lambda}b_{2}))+\delta(b_{1}\oo{\lambda}b_{2})
+\hd_{A}(\chi_{\lambda}(b_{1}, b_{2}))\\
=&\; \gamma(\hd_{B}(b_{1})\oo{\lambda}b_{2})+\gamma(b_{1}\oo{\lambda}\hd_{B}(b_{2}))
+\delta(b_{1}\oo{\lambda}b_{2})+\hd_{A}(\chi_{\lambda}(b_{1}, b_{2})),
\end{align*}
and
\begin{align*}
&\; \hd_{E}(\gamma(b_{1}))\oo{\lambda}\gamma(b_{2})
+\gamma(b_{1})\oo{\lambda}\hd_{E}(\gamma(b_{2}))\\
=&\; \Big(\gamma(\hd_{B}(b_{1}))+\delta(b_{1})\Big)\oo{\lambda}\gamma(b_{2})
+\gamma(b_{1})\oo{\lambda}\Big(\gamma(\hd_{B}(b_{2}))+\delta(b_{2})\Big)\\
=&\; \gamma(\hd_{B}(b_{1}))\oo{\lambda}\gamma(b_{2})+\delta(b_{1})\oo{\lambda}\gamma(b_{2})
+\gamma(b_{1})\oo{\lambda}\gamma(\hd_{B}(b_{2}))+\gamma(b_{1})\oo{\lambda}\delta(b_{2}).
\end{align*}
Since $\hd_{E}$ is a derivation, i.e., $\hd_{E}\Big(\gamma(b_{1})\oo{\lambda}\gamma(b_{2})\Big)
=\hd_{E}(\gamma(b_{1}))\oo{\lambda}\gamma(b_{2})+\gamma(b_{1})\oo{\lambda}\hd_{E}(\gamma(b_{2}))$,
we get
\begin{align*}
&\; \hd_{E}(\chi_{\lambda}(b_{1}, b_{2}))-\chi_{\lambda}(\hd_{B}(b_{1}), b_{2})
-\chi_{\lambda}(b_{1}, \hd_{B}(b_{2}))\\
=&\; \gamma(b_{1})\oo{\lambda}\delta(b_{2})-\delta(b_{1}\oo{\lambda}b_{2})
+\delta(b_{1})\oo{\lambda}\gamma(b_{2}).
\end{align*}
That is to say, for any $b_{1}, b_{2}\in B$,
\begin{align*}
\Theta(\hd_{A}, \hd_{B})(\chi)_{\lambda}(b_{1}, b_{2})
&=\hd_{E}(\chi_{\lambda}(b_{1}, b_{2}))-\chi_{\lambda}(\hd_{B}(b_{1}), b_{2})
-\chi_{\lambda}(b_{1}, \hd_{B}(b_{2}))\\
&=d_{1}(\delta)_{\lambda}(b_{1}, b_{2}).
\end{align*}
Hence $\mathcal{W}(\hd_{A}, \hd_{B})=[0]$. The proof is finished.
\end{proof}

Now we consider the exact sequence of Wells about derivations. Let $B$
be an associative conformal algebra, $A$ be a $B$-bimodule, and
$\mathcal{E}: \xymatrix@C=0.5cm{0\ar[r]& A\ar[r]^{\alpha}&E\ar[r]^{\beta}&B\ar[r]&0}$
be an abelian extension of $B$ by $A$ with a section $\gamma$. Denote
$$
\Der_{A}(E):=\{\hd_{E}\in\Der(E)\mid \hd_{E}(A)\subseteq A\}.
$$
Then $\Der_{A}(E)$ is a Lie subalgebra of $\Der(E)$. Any $\hd_{E}\in\Der_{A}(E)$,
induces two maps $\hd_{E}|_{A}\in\Der(A)$ and $\hat{\hd}_{E}=\beta\circ\hd_{E}\circ\gamma:
B\rightarrow B$. Then one can check that $\hat{\hd}_{E}\in\Der(B)$ by direct calculation
and it does not depend on the choice of $\gamma$. Thus we get a linear map
$$
\kappa:\; \Der_{A}(E)\rightarrow\Der(A)\times\Der(B),\qquad\qquad
\hd_{E}\mapsto(\hd_{E}|_{A}, \hat{\hd}_{E}).
$$

\begin{lem} \label{lem: liemap}
With the above notations, $\Img(\kappa)\subseteq\g(A, B)$, and $\kappa$ induces
a Lie algebra homomorphism $\bar{\kappa}: \Der_{A}(E)\rightarrow\g(A, B)$.
\end{lem}

\begin{proof}
Let $(\hd_{A}, \hd_{B})\in\Img(\kappa)$. There exists a derivation $\hd_{E}\in\Der(E)$
such that $\hd_{A}=\hd_{E}|_{A}$ and $\hd_{B}=\hat{\hd}_{E}$. Then clearly,
$\hd_{E}\circ\alpha=\alpha\circ\hd_{A}$. Moreover, for any $a\in A$ and $b\in B$, we have
$$
\hat{\hd}_{E}(\beta(a+\gamma(b)))=\hat{\hd}_{E}(b)=\beta(\hd_{E}(a+\gamma(b))).
$$
That is, $\hd_{B}\circ\beta=\beta\circ\hd_{E}$. Hence $\Img(\kappa)\subseteq\g(A, B)$,
and so that $\kappa$ induces a linear map $\bar{\kappa}: \Der_{A}(E)\rightarrow\g(A, B)$.
We now show that $\bar{\kappa}$ is a homomorphism of Lie algebras.
Let $\hd_{E}$, $\hd'_{E}\in\Der_{A}(E)$. For any $b\in B$, there exist $a_{1}, a_{2}\in A$
and $b_{1}, b_{2}\in B$ such that
$$
\hd_{E}(\gamma(b))=a_{1}+\gamma(b_{1}),\qquad\qquad \hd'_{E}(\gamma(b))=a_{2}+\gamma(b_{2}).
$$
Then $\hat{\hd}_{E}(b)=b_{1}$, $\hat{\hd}'_{E}(b)=b_{2}$, and hence
\begin{align*}
\widehat{[\hd_{E}, \hd'_{E}]}(b)
&=\beta\Big((\hd_{E}\circ\hd'_{E}-\hd'_{E}\circ\hd_{E})(\gamma(b))\Big)\\
&=\beta\Big(\hd_{E}(a_{2}+\gamma(b_{2}))\Big)+\beta\Big(\hd'_{E}(a_{1}+\gamma(b_{1}))\Big)\\
&=\hat{\hd}_{E}(b_{2})-\hat{\hd}'_{E}(b_{1})\\
&=[\hat{\hd}_{E},\; \hat{\hd}'_{E}](b).
\end{align*}
Therefore,
\begin{align*}
\bar{\kappa}([\hd_{E}, \hd'_{E}])
&=\Big([\hd_{E}, \hd'_{E}]|_{A},\; \widehat{[\hd_{E}, \hd'_{E}]}\Big)\\
&=\Big([\hd_{E}|_{A}, \hd'_{E}|_{A}],\; [\hat{\hd}_{E},\; \hat{\hd}'_{E}]\Big)\\
&=[(\hd_{E}|_{A}, \hat{\hd}_{E}),\; (\hd'_{E}|_{A}, \hat{\hd}'_{E})]\\
&=[\bar{\kappa}(\hd_{E}), \; \bar{\kappa}(\hd'_{E})].
\end{align*}
Thus $\bar{\kappa}$ is a homomorphism of Lie algebras.
\end{proof}

Consider the image and kernel of $\bar{\kappa}$, we get an exact sequence of Wells type.

\begin{thm}\label{thm: well-der}
Let $\mathcal{E}: \xymatrix@C=0.5cm{0\ar[r]& A\ar[r]^{\alpha}&E\ar[r]^{\beta}&B\ar[r]&0}$
be an abelian extension of associative conformal algebra $B$ by bimodule $A$.
There is an exact sequence of vector spaces
$$
\xymatrix{0\ar[r]&\mathcal{Z}^{1}(B, A)\ar[r]^{\iota}&\Der_{A}(E)\ar[r]^{\bar{\kappa}}
&\g(A, B)\ar[r]^{\mathcal{W}\ \ } &\HH^{2}(B, A)},
$$
where $\iota$ is the inclusion map.
\end{thm}

\begin{proof}
If $\hd_{E}\in\Ker(\bar{\kappa})$, then $\hat{\hd}_{E}=0$ and $\hd_{E}|_{A}=0$.
For any $b\in B$, we get $\beta(\hd_{E}(\gamma(b)))=\hat{\hd}_{E}(b)=0$, that is,
$\hd_{E}(\gamma(b))\in A$. Thus there is a linear map
$$
\phi:\; \Ker(\bar{\kappa})\longrightarrow\mathcal{C}^{1}(B, A)=\Hom_{\cb[\partial]}(B, A),
\qquad\qquad \hd_{E}\mapsto \hd_{E}\circ\gamma.
$$
Then one can check that $\phi$ does not depend on the choice of $\gamma$, and
$\phi(\hd_{E})\in\mathcal{Z}^{1}(B, A)$ for any $\hd_{E}\in\Ker(\bar{\kappa})$.
Thus $\phi$ induces a linear map $\bar{\phi}: \Ker(\bar{\kappa})\rightarrow
\mathcal{Z}^{1}(B, A)$. First, if $\hd_{E}\in \Ker(\bar{\kappa})$ and $\bar{\phi}(\hd_{E})=0$,
for any $a\in A$ and $b\in B$, we get
$$
\hd_{E}(a+\gamma(b))=\hd_{E}|_{A}(a)+\bar{\phi}(\hd_{E})(b)=0.
$$
This means that $\bar{\phi}$ is injective. Second, for any $f\in\mathcal{Z}^{1}(B, A)$,
we define a $\cb[\partial]$-module homomorphism
$$
\hd_{E}:\; E\longrightarrow E,\qquad\qquad    \hd_{E}(a+\gamma(b))=f(b),
$$
for any $a\in A$ and $b\in B$. Clearly, $\hd_{E}|_{A}=0$, $\hat{\hd}_{E}(b)=\beta(f(b))=0$
and $\bar{\phi}(\hd_{E})(b)=\hd_{E}(\gamma(b))=f(b)$, for any $b\in B$.
Moreover, one can check that $\hd_{E}\in\Der(E)$, and so that $\hd_{E}\in\Ker(\bar{\kappa})$.
Thus $\bar{\phi}$ is surjective, and is an isomorphism of vector spaces.

Next, we show that the sequence is exact at $\g(A, B)$. For any $(\hd_{A}, \hd_{B})\in
\Img(\bar{\kappa})$, which is extensible, we get $\mathcal{W}(\hd_{A}, \hd_{B})=[0]$
by Theorem \ref{thm: ext-con}. Thus $\Img(\bar{\kappa})\subseteq\Ker(\mathcal{W})$.
Conversely, for any $(\hd_{A}, \hd_{B})\in\Ker(\mathcal{W})$, that is, there exists a
element $f\in\mathcal{C}^{1}(B, A)$ such that $\Theta(\hd_{A}, \hd_{B})([\chi])=
[d_{1}(f)]$. Similar to the proof of Theorem \ref{thm: ext-con}, we define
$\hd_{E}: E\rightarrow E$ by
$$
\hd_{E}(a+\gamma(b))=\hd_{A}(a)+f(b)+\gamma(\hd_{B}(b)),
$$
for any $a\in A$ and $b\in B$. Then $\hd_{E}\in\Der(E)$ and $\bar{\kappa}(\hd_{E})
=(\hd_{A}, \hd_{B})$. That is, $\Ker(\mathcal{W})\subseteq\Img(\bar{\kappa})$.
Therefore, we get the sequence which is an exact sequence.
\end{proof}

Finally, we consider a special case, that is, when the abelian extension of associative
conformal algebras is split. In this case, we can see that there is a close relationship
between the Lie algebra $\Der_{A}(E)$ and $\g(A, B)$.

\begin{cor}\label{cor:split-der}
Let $\mathcal{E}: \xymatrix@C=0.5cm{0\ar[r]& A\ar[r]^{\alpha}&E\ar[r]^{\beta}&B\ar[r]&0}$
be a split abelian extension of associative conformal algebra $B$ by bimodule $A$, i.e.,
there exists a homomorphism of associative conformal algebras $\gamma: B\rightarrow E$
such that $\beta\circ\gamma=\id_{B}$. Then
$$
\Der_{A}(E)\cong \g(A, B)\times\mathcal{Z}^{1}(B, A)
$$
as Lie algebras, where $\mathcal{Z}^{1}(B, A)$ is regarded as a trivial Lie algebra.
\end{cor}

\begin{proof}
Since $\mathcal{E}$ is split, by Theorem \ref{thm: well-der}, we get a short exact
sequence of Lie algebras
$$
(\star)\qquad\qquad\xymatrix{0\ar[r]&\mathcal{Z}^{1}(B, A)\ar[r]^{\iota}&\Der_{A}(E)
\ar[r]^{\bar{\kappa}}&\g(A, B)\ar[r]^{\; \mathcal{W} } &0},\qquad\qquad\quad
$$
if we regard $\mathcal{Z}^{1}(B, A)$ as a trivial Lie algebra. For any
$(\hd_{A}, \hd_{B})\in\g(A, B)$, we define
$$
\hd_{E}:\; E\longrightarrow E,\qquad\qquad
\hd_{E}(a+\gamma(b))=\hd_{A}(a)+\gamma(\hd_{B}(b)),
$$
for any $a\in A$ and $b\in B$. Note that in this case the cohomology class $[\chi]=0$,
by the proof of Theorem \ref{thm: ext-con}, we get $\hd_{E}\in\Der(E)$. Thus we obtain a
linear map
$$
\eta:\; \g(A, B)\longrightarrow\Der(E),\qquad\qquad
(\hd_{A}, \hd_{B})\mapsto\hd_{E}.
$$
It is easy to see $\bar{\kappa}\circ\eta(\hd_{A}, \hd_{B})=\bar{\kappa}(\hd_{E})
=(\hd_{E}|_{A},\; \hat{\hd}_{E})=(\hd_{A}, \hd_{B})$. That is,
$\bar{\kappa}\circ\eta=\id_{\g(A, B)}$. Moreover, by direct calculation, one can show
that $\eta$ is a homomorphism of Lie algebras. Thus the sequence $(\star)$ is split
in the category of Lie algebras. Hence we get $\Der_{A}(E)\cong \g(A, B)\times
\mathcal{Z}^{1}(B, A)$ as Lie algebras.
\end{proof}

\begin{rmk}
In the context of associative conformal algebras, there is a class of derivation,
called conformal derivation (see \cite{BKV}). All conformal derivations of an
associative conformal algebra form a Lie conformal algebra.
But the cohomology of associative conformal algebra we use here has only a vector space
structure, no $\cb[\partial]$-module structure. This cohomology cannot form a
representation of a Lie conformal algebra, and no relevant conclusions can be obtained.
Therefore, here we still consider derivations on associative conformal algebras
instead of conformal derivations.
\end{rmk}

\bigskip
\noindent
{\bf Acknowledgements. } This work was financially supported by National
Natural Science Foundation of China (No. 11771122, 11801141, 12201182), China
Postdoctoral Science Foundation (2020M682272) and NSF of Henan Province
(212300410120).

 \end{document}